\documentclass[11pt]{article}

\usepackage[utf8]{inputenc} %
\usepackage[T1]{fontenc}    %
\usepackage{hyperref}       %
\usepackage{url}            %
\usepackage{booktabs}       %
\usepackage{amsfonts}       %
\usepackage{nicefrac}       %
\usepackage{microtype}      %
\usepackage{amsmath}
\usepackage{bm}
\usepackage{amssymb}
\usepackage{mathtools}
\usepackage{graphicx}
\usepackage{physics}
\usepackage{float}
\usepackage{yhmath}
\usepackage{algorithm,algorithmic}

\usepackage{tikz}

\usepackage{amsthm}

\usepackage{empheq}

\usepackage{natbib}

\let\savedegree\div
\let\div\relax
\let\savering\ring
\let\ring\relax
\usepackage{mathabx} %
\let\div\savedegree
\let\ring\savering

\newcommand{\E}[1]{\mathbb{E}\left[#1\right]}

\newcommand{\bracket}[1]{\left[ #1 \right]}

\newcommand{\bg}[1]{\boldsymbol{#1}}
\newcommand{\balpha}{\boldsymbol{\alpha}}
\newcommand{\btheta}{\boldsymbol{\theta}}
\newcommand{\bW}{\boldsymbol{W}}

\newcommand{\EE}{\mathbb{E}}
\newcommand{\RR}{\mathbb{R}}

\newcommand{\mc}[1]{\mathcal{#1}}

\newtheorem{remark}{Remark}

\theoremstyle{definition}
\newtheorem{exmp}{Example}[section]

\DeclareMathOperator*{\argmax}{arg\,max}
\DeclareMathOperator*{\argmin}{arg\,min}

\title{Reinforcement Learning for Mean Field Games, with Applications to Economics}

\author{
Andrea Angiuli\thanks{Department of Statistics and Applied Probability,
South Hall 5504,
University of California
Santa Barbara, CA 93106
{\tt angiuli@pstat.ucsb.edu}} 
\and 
Jean-Pierre Fouque\thanks{
Department of Statistics and Applied Probability,
South Hall 5504,
University of California
Santa Barbara, CA 93106
{\tt fouque@pstat.ucsb.edu}. Work  supported by NSF grant DMS-1814091.
}
\and Mathieu Lauri{\`e}re\thanks{Department of Operations Research and Financial Engineering. Princeton University. {\tt lauriere@princeton.edu}. Work supported by ARO grant AWD1005491 and NSF award AWD1005433.}}

\begin{document}
\maketitle

\begin{abstract}
Mean field games (MFG) and mean field control problems (MFC) are frameworks to study Nash equilibria or social optima in games with a continuum of agents. These problems can be used to approximate competitive or cooperative games with a large finite number of agents and have found a broad range of applications, in particular in economics. In recent years, the question of learning in MFG and MFC has garnered interest, both as a way to compute solutions and as a way to model how large populations of learners converge to an equilibrium. Of particular interest is the setting where the agents do not know the model, which leads to the development of reinforcement learning (RL) methods. After reviewing the literature on this topic, we present a two timescale approach with RL for MFG and MFC, which relies on a unified Q-learning algorithm. The main novelty of this method is to simultaneously update an action-value function and a distribution but with different rates, in a model-free fashion. Depending on the ratio of the two learning rates, the algorithm learns either the MFG or the MFC solution. To illustrate this method, we apply it 
to a mean field problem of accumulated consumption in finite horizon with HARA utility function, and to 
a trader's optimal liquidation problem. 
\end{abstract}

\section{Introduction}

Dynamic games with many players are pervasive in today's highly connected world. In many models the agents are indistinguishable since they have the same dynamics and cost functions. Moreover, the interactions are often anonymous since each player is influenced only by the empirical distribution of all the agents. However, such games become intractable when the number of agents becomes very large. Mean field games have been introduced in~\citep{MR2295621,MR2346927} to tackle such situations by passing to the limit and considering games with an infinite number of players interacting through the population distribution. Although the standard formulation of MFG focuses on finding Nash equilibria, social optima arising in a cooperative setting have also been studied under the term of mean field control~\citep{bensoussan2013mean} or control of McKean-Vlasov dynamics~\citep{MR2295621}. Equilibria or social optima in such games can be characterized in a tractable way through forward-backward systems of partial differential equations (PDE) or stochastic differential equations (SDE)~\citep{MR3752669,MR2295621}.

Mean field games with interactions through the controls, sometimes called ``extended'', occur when the dynamics or the cost function of a typical player explicitly depends on the empirical measure of the \emph{controls} of the other players, and not just on their respective states.
Such games were first introduced by Gomes et. al. \citep{MR3160525,MR3489048} and their investigation quickly garnered interest.

Interaction through the controls' distribution is particularly relevant in economics and finance, see \textit{e.g.}~\citep{huang2013mean,MR3195844,MR3325272,MR3359708,MR3755719,MR3805247} and ~\citep{carmona-AMS-tmp} for a recent survey.  
Some aspects of the PDE approach and the probabilistic approach to such games have been treated respectively in~\citep{MR3941633,BonnansHadikhanlooPfeiffer2019,Kobeissi2019} and in~\citep{MR3325272}. 
As in many fields, linear-quadratic models are particularly appealing due to their tractability, see \textit{e.g.}~\citep{alasseur2017extended,MR3575612} for applications to energy production. 

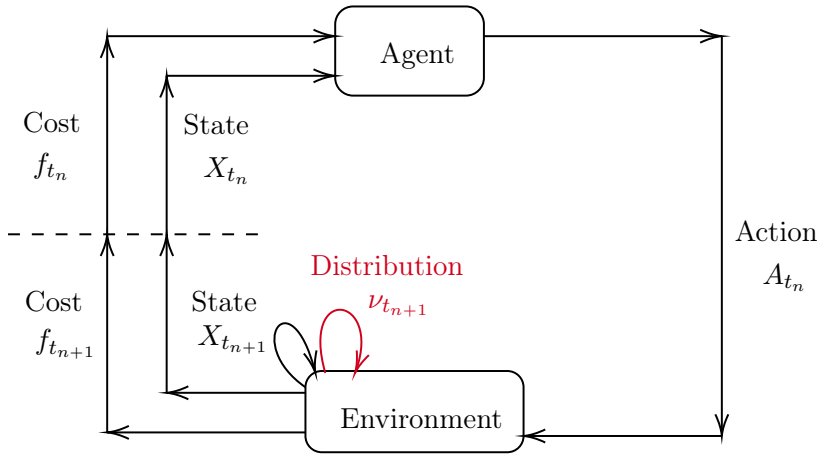
\begin{figure}
\begin{center}

\tikzset{every picture/.style={line width=0.75pt}} %

\begin{tikzpicture}[x=0.75pt,y=0.75pt,yscale=-1,xscale=1]

\draw   (220,227.2) .. controls (220,222.67) and (223.67,219) .. (228.2,219) -- (321.8,219) .. controls (326.33,219) and (330,222.67) .. (330,227.2) -- (330,251.8) .. controls (330,256.33) and (326.33,260) .. (321.8,260) -- (228.2,260) .. controls (223.67,260) and (220,256.33) .. (220,251.8) -- cycle ;
\draw    (220,250) -- (122,250) ;
\draw [shift={(120,250)}, rotate = 360] [color={rgb, 255:red, 0; green, 0; blue, 0 }  ][line width=0.75]    (10.93,-3.29) .. controls (6.95,-1.4) and (3.31,-0.3) .. (0,0) .. controls (3.31,0.3) and (6.95,1.4) .. (10.93,3.29)   ;
\draw    (120,250) -- (120,152) ;
\draw [shift={(120,150)}, rotate = 450] [color={rgb, 255:red, 0; green, 0; blue, 0 }  ][line width=0.75]    (10.93,-3.29) .. controls (6.95,-1.4) and (3.31,-0.3) .. (0,0) .. controls (3.31,0.3) and (6.95,1.4) .. (10.93,3.29)   ;
\draw   (235,44) .. controls (235,39.03) and (239.03,35) .. (244,35) -- (301,35) .. controls (305.97,35) and (310,39.03) .. (310,44) -- (310,71) .. controls (310,75.97) and (305.97,80) .. (301,80) -- (244,80) .. controls (239.03,80) and (235,75.97) .. (235,71) -- cycle ;
\draw    (120,50) -- (233,50) ;
\draw [shift={(235,50)}, rotate = 180] [color={rgb, 255:red, 0; green, 0; blue, 0 }  ][line width=0.75]    (10.93,-3.29) .. controls (6.95,-1.4) and (3.31,-0.3) .. (0,0) .. controls (3.31,0.3) and (6.95,1.4) .. (10.93,3.29)   ;
\draw    (220,230) -- (152,230) ;
\draw [shift={(150,230)}, rotate = 360] [color={rgb, 255:red, 0; green, 0; blue, 0 }  ][line width=0.75]    (10.93,-3.29) .. controls (6.95,-1.4) and (3.31,-0.3) .. (0,0) .. controls (3.31,0.3) and (6.95,1.4) .. (10.93,3.29)   ;
\draw    (150,230) -- (150,152) ;
\draw [shift={(150,150)}, rotate = 450] [color={rgb, 255:red, 0; green, 0; blue, 0 }  ][line width=0.75]    (10.93,-3.29) .. controls (6.95,-1.4) and (3.31,-0.3) .. (0,0) .. controls (3.31,0.3) and (6.95,1.4) .. (10.93,3.29)   ;
\draw    (150,70) -- (233,70) ;
\draw [shift={(235,70)}, rotate = 180] [color={rgb, 255:red, 0; green, 0; blue, 0 }  ][line width=0.75]    (10.93,-3.29) .. controls (6.95,-1.4) and (3.31,-0.3) .. (0,0) .. controls (3.31,0.3) and (6.95,1.4) .. (10.93,3.29)   ;
\draw    (430,249.8) -- (430,50) ;
\draw [shift={(430,251.8)}, rotate = 270] [color={rgb, 255:red, 0; green, 0; blue, 0 }  ][line width=0.75]    (10.93,-3.29) .. controls (6.95,-1.4) and (3.31,-0.3) .. (0,0) .. controls (3.31,0.3) and (6.95,1.4) .. (10.93,3.29)   ;
\draw    (310,50) -- (428,50) ;
\draw [shift={(430,50)}, rotate = 180] [color={rgb, 255:red, 0; green, 0; blue, 0 }  ][line width=0.75]    (10.93,-3.29) .. controls (6.95,-1.4) and (3.31,-0.3) .. (0,0) .. controls (3.31,0.3) and (6.95,1.4) .. (10.93,3.29)   ;
\draw    (430,251.8) -- (332,251.8) ;
\draw [shift={(330,251.8)}, rotate = 360] [color={rgb, 255:red, 0; green, 0; blue, 0 }  ][line width=0.75]    (10.93,-3.29) .. controls (6.95,-1.4) and (3.31,-0.3) .. (0,0) .. controls (3.31,0.3) and (6.95,1.4) .. (10.93,3.29)   ;
\draw  [dash pattern={on 4.5pt off 4.5pt}]  (70,150) -- (200,150) ;
\draw    (220,227.2) .. controls (188.32,205.22) and (212.51,171.68) .. (224.64,218.55) ;
\draw [shift={(225,220)}, rotate = 256.24] [color={rgb, 255:red, 0; green, 0; blue, 0 }  ][line width=0.75]    (10.93,-3.29) .. controls (6.95,-1.4) and (3.31,-0.3) .. (0,0) .. controls (3.31,0.3) and (6.95,1.4) .. (10.93,3.29)   ;
\draw    (150,150) -- (150,72) ;
\draw [shift={(150,70)}, rotate = 450] [color={rgb, 255:red, 0; green, 0; blue, 0 }  ][line width=0.75]    (10.93,-3.29) .. controls (6.95,-1.4) and (3.31,-0.3) .. (0,0) .. controls (3.31,0.3) and (6.95,1.4) .. (10.93,3.29)   ;
\draw    (120,150) -- (120,52) ;
\draw [shift={(120,50)}, rotate = 450] [color={rgb, 255:red, 0; green, 0; blue, 0 }  ][line width=0.75]    (10.93,-3.29) .. controls (6.95,-1.4) and (3.31,-0.3) .. (0,0) .. controls (3.31,0.3) and (6.95,1.4) .. (10.93,3.29)   ;
\draw [color={rgb, 255:red, 208; green, 2; blue, 27 }  ,draw opacity=1 ][line width=0.75]    (230,220) .. controls (217.2,174.69) and (258.72,180.8) .. (245.64,218.27) ;
\draw [shift={(245,220)}, rotate = 291.04] [color={rgb, 255:red, 208; green, 2; blue, 27 }  ,draw opacity=1 ][line width=0.75]    (10.93,-3.29) .. controls (6.95,-1.4) and (3.31,-0.3) .. (0,0) .. controls (3.31,0.3) and (6.95,1.4) .. (10.93,3.29)   ;

\draw (236,236) node [anchor=north west][inner sep=0.75pt]   [align=left] {Environment};
\draw (256,52) node [anchor=north west][inner sep=0.75pt]   [align=left] {Agent};
\draw (77,177) node [anchor=north west][inner sep=0.75pt]   [align=left] {Cost};
\draw (82,196.4) node [anchor=north west][inner sep=0.75pt]    {$f_{t_{n+1}}$};
\draw (161,177) node [anchor=north west][inner sep=0.75pt]   [align=left] {State};
\draw (164,194.4) node [anchor=north west][inner sep=0.75pt]    {$X_{t_{n+1}}$};
\draw (221,159) node [anchor=north west][inner sep=0.75pt]  [color={rgb, 255:red, 208; green, 2; blue, 27 }  ,opacity=1 ] [align=left] {Distribution};
\draw (250,179.4) node [anchor=north west][inner sep=0.75pt]  [color={rgb, 255:red, 208; green, 2; blue, 27 }  ,opacity=1 ]  {$\nu _{t_{n+1}}$};
\draw (435,142) node [anchor=north west][inner sep=0.75pt]   [align=left] {Action};
\draw (449,163.4) node [anchor=north west][inner sep=0.75pt]    {$A_{t_{n}}$};
\draw (76,87) node [anchor=north west][inner sep=0.75pt]   [align=left] {Cost};
\draw (81,107.4) node [anchor=north west][inner sep=0.75pt]    {$f_{t_{n}}$};
\draw (157,88) node [anchor=north west][inner sep=0.75pt]   [align=left] {State};
\draw (167,109.4) node [anchor=north west][inner sep=0.75pt]    {$X_{t_{n}}$};

\end{tikzpicture}

\end{center}

\caption{At time $t_n$ of an experiment, the agent feeds the action $A_{t_n}$ to the environment, which outputs the associated cost $f_{t_{n+1}} = f(X_{t_n},A_{t_n}, \nu_{t_n})$ for $n < N_T$ (and $f_{t_{n+1}} = g(X_{t_n}, \nu_{t_n})$ for $n = N_T$) and the new state $X_{t_{n+1}} \sim p(X_{t_n}, A_{t_n},\nu_{t_n})$. Besides these outputs provided to the agent, the environment keeps track of $X_{t_{n+1}}$ and $\nu_{t_{n+1}}$ for the next iteration. }
\label{diagram}
\end{figure}

The approach we propose is based on ideas from reinforcement learning (RL). Applications of RL in economics and finance have recently attracted a lot of interest, see \textit{e.g.}~\citep{charpentier2020reinforcement}. However, since our problems involve mean-field interactions, the population distribution requires a special treatment. 
In our setup, the agent is feeding an action to the environment
which produces the next state and a reward (or cost). The environment also updates in an automatic way (without decision) the distributions of states and controls. Then, the agent update their Q-matrix and proceeds (see the diagram in Figure \ref{diagram}). The environment can be viewed as a ``black box" or as a ``simulator'' depending on the problem, but, in any case, it generates the new state if the dynamics is unknown and the reward if not computable by the agent. It is also interesting to note that even in cases where the dynamics and the reward structure are known but complicated, then our algorithm can be viewed as a numerical method for computing the optimal strategy for the corresponding MFG or MFC problems.

Since the introduction of MFG theory, several numerical methods have been proposed, see \textit{e.g.}~\citep{achdoulauriere-2020-mfg-numerical,AMSnotesLauriere} and the references therein. Recently, several methods to solve MFGs based on machine learning tools have been proposed relying either on the probabilistic approach~\citep{fouque2019deep,carmona2019convergence-II,germain2019numerical,min2021signatured} or the analytical approach~\citep{al2018solving,carmona2021convergence,ruthotto2020machine,cao2020connecting,lin2020apac,AMSnotesLauriere}. They combine neural network approximations and stochastic optimization techniques to solve McKean-Vlasov control problems, mean field FBSDE or mean field PDE systems; see ~\citep{CarmonaLauriere-2021-DLMFGFIN} for a recent survey and applications to finance. These methods are based on the knowledge of the model, but the question of learning solutions to MFG and MFC without full knowledge of the model have also attracted a surge of interest. 

As far as learning methods for mean field problems are concerned, most works focus either on MFG or on MFC. \citep{yang2018mean} uses a mean field approximation in the context of multi-agent reinforcement learning (MARL) to reduce the computational cost. \citep{yang2018deep} uses inverse reinforcement learning to learn the dynamics of a mean field game on a graph. To approximate stationary MFG solutions, \citep{guo2019learning} uses fixed point iterations on the distribution combined with Q-learning to learn the best response at each iteration. \citep{anahtarci2020q} combines this kind of learning scheme together with an entropic regularization. Convergence of an actor-critic method for linear-quadratic MFG has been studied in~\citep{fu2019actor}. Model-free learning for finite horizon MFG has been studied in~\citep{mishra2020model} using a backward scheme. Fictitious play without or with reinforcement learning has been studied respectively in~\citep{cardaliaguet2017learning,hadikhanloo2019finite} and~\citep{elie2020convergence,perrin2020continuousfp,xie2020provable}, or online mirror descent~\citep{hadikhanloo2017learning,perolat2021OMD}. These iterative methods have been proved to converge under a monotonicity condition which is weaker than the strict contraction property used to ensure convergence of fixed point iterations. They can be extended to continuous space problems using deep reinforcement learning as \textit{e.g.} in~\citep{perrin2021mean}. A two timescale approach to solve MFG with finite state and action spaces has been proposed in~\citep{mguni2018decentralised,SubramanianMahajan-2018-RLstatioMFG}. 

To learn MFC optima, \citep{SubramanianMahajan-2018-RLstatioMFG} designs a gradient based algorithm. Model-free policy gradient method has been proved to converge for linear-quadratic problems in~\citep{CarmonaLauriereTan-2019-LQMFRL,wang2020global}, whereas Q-learning for a ``lifted'' Markov decision process on the space of distributions has been studied  in~\citep{CarmonaLauriereTan-2019-MQFRL,gu2019dynamic,gu2020q}. Optimality conditions and propagation of chaos type result for mean field Markov decision processes  are studied in~\citep{motte2019mean}.

In~\citep{angiuli2020unified}, we proposed a unified two timescale Q-learning algorithm to solve both MFG and MFC problems in an infinite horizon stationary regime. The key idea is to iteratively update estimates of the distribution and the Q-function with different learning rates. Suitably choosing these learning rates enables the algorithm to learn the solution of the MFG or the one of the MFC. A slow updating of the distribution of the state leads to the Nash equilibrium of the competitive MFG and the algorithm learns the corresponding optimal strategy. A rapid updating of the distribution leads to learning of the optimal control  of the corresponding cooperative MFC. Moreover, in contrast with other approaches, our algorithm does not require the environment to output the population distribution which means that a single agent can learn the solution of mean field problems.

In the present work, we extend this algorithm in two directions: finite horizon setting, and ``extended" mean field problems which involve the distribution of controls as well. That demonstrates the flexibility of our two timescale algorithm and broadens the range of applications.

The rest of the paper is organized as follows. In Section~\ref{sec:FTMFGMFC}, we introduce the framework of finite horizon mean field games and mean field control problems. In Section~\ref{sec:twotimescale}, we present the main ideas behind the two timescale approach in this context. Based on this perspective, we introduce in Section~\ref{sec:RLalgo} a reinforcement learning algorithm to solve MFC and MFG problems. We then illustrate this method on two examples: a mean field accumulation problem in Section~\ref{sec:accumulation} and an optimal execution problem for a mean field of traders in Section~\ref{sec:trader}. We then conclude in Section~\ref{sec:conclusion}.

\noindent{\bf Notation. } For a random variable $X$, $\mathcal{L}(X)$ denotes its law. $d$ and $k$ are two positive integers corresponding respectively to the state and the action dimensions. Unless otherwise specified, $\nu$ will be used to denote a state-action distribution, and its first and second marginals will respectively be denoted by $\mu$ and $\theta$.

\section{Finite horizon mean field problems}\label{sec:FTMFGMFC}

In this section, we introduce the framework of mean field games and mean field control problems in finite horizon. For the sake of consistency with the MFG literature, we use a continuous time formalism. For the link with finite player games, see \textit{e.g.}~\citep{MR3752669}.

\subsection{Mean field games} \label{sec:mfg}

Let $(\Omega,\mathcal{F},\mathbb{F}=(\mathcal{F}_t)_{0 \leq t \leq T},\mathbb{P})$ be a filtered probability space, where the filtration supports an $m$-dimensional Brownian motion $W=(W_t)_{0\leq t \leq T}$ and an initial condition $\xi \in L^2(\Omega,\mathcal{F}_0,\mathbb{P};\mathbb{R}^d)$. Let $f: [0,T] \times \mathbb{R}^d \times \mathcal{P}_2(\mathbb{R}^{d+k}) \times \mathbb{R}^k \to \mathbb{R}$ and $g: \mathbb{R}^d \times \mathcal{P}_2(\mathbb{R}^{d}) \to \mathbb{R}$ be respectively a running cost function and a terminal cost function. Let $b: [0,T] \times \mathbb{R}^d \times \mathcal{P}_2(\mathbb{R}^{d+k}) \times \mathbb{R}^k \to \mathbb{R}^d$ be a drift function and let $\sigma: [0,T] \times \mathbb{R}^d \times \mathcal{P}_2(\mathbb{R}^{d+k}) \times \mathbb{R}^k \to \mathbb{R}^{d \times m}$ be a volatility function.

A mean field game equilibrium is defined as a pair $$(\hat\alpha,\hat\nu) = ((\hat\alpha_t)_{t \in [0,T]},(\hat\nu_t)_{t \in [0,T]}) \in \mathbb{A} \times \mathcal{C}([0,T] , \mathcal{P}_2(\mathbb{R}^{d+k})),$$ where $\mathbb{A}$ is the set of admissible controls, namely progressively measurable processes that are square integrable, such that

\begin{enumerate}
    \item $\hat\alpha$ solve the standard stochastic control problem when $\nu = \hat\nu$:
    \begin{equation*}
    \inf_{\alpha \in \mathbb{A}}J_{\nu}(\alpha)=\inf_{\alpha \in \mathbb{A}}\mathbb{E}\left[\int_0^T f(t,X^{\alpha,\nu}_t,\nu_t,\alpha_t)dt+g(X^{\alpha,\nu}_T,\hat\mu_T) \right],
    \end{equation*}
    where $\hat\mu_T$ is the first marginal of $\nu_T$ corresponding to the terminal state distribution 
    subject to
    \begin{equation*}
    dX^{\alpha,\nu}_t=b(t,X^{\alpha,\nu}_t,\nu_t,\alpha_t)dt+\sigma(t,X^{\alpha,\nu}_t,\nu_t,\alpha_t) dW_t,\quad
    X^{\alpha,\nu}_0=\xi.
    \end{equation*}
    
    \item $\hat\nu_t = \mathcal{L}(X^{\hat\alpha}_t, \hat\alpha_t)$ for all $0\leq t \leq T$. 
\end{enumerate}

The solution can be characterized from either the PDE viewpoint (leading to a coupled Hamilton-Jacobi-Bellman and Kolmogorov-Fokker-Plank equations) \citep{MR2295621,MR2346927} or from a probabilistic viewpoint \citep{carmona2018probabilisticI-II}. Within the probabilistic viewpoint, there are two approaches, both of which are formulated with FBSDEs. The backward variable can represent either the value function of a typical player or the derivative of this value function. See~\cite[Volume 1, Chapters 3 and 4]{carmona2018probabilisticI-II} for more details. These analytical and probabilistic approaches also lead to computational methods, as long as the model is known. However, when the model is not known, one needs to develop other tools, as we will discuss in the next sections.

\subsection{Mean field control}\label{sec:mfc}

In contrast with the MFG problem -- which corresponds to a Nash equilibrium, the mean field control (MFC) problem is an optimization problem. It can be interpreted as the problem posed to a social planner trying to find the optimal behavior of a population so as to minimize a social cost (\textit{i.e.}, a cost averaged over the whole population). It is an optimal control problem for a McKean-Vlasov dynamics: Find $\alpha^*$ which satisfies:
    \begin{equation*}
    \inf_{\alpha \in \mathbb{A}}J(\alpha)=\inf_{\alpha \in \mathbb{A}}\mathbb{E}\left[\int_0^T f(t,X^{\alpha}_t,\nu^{\alpha}_t,\alpha_t)dt+g(X^{\alpha}_T,\mu^{\alpha}_T) \right],
    \end{equation*}
    subject to
    
    \begin{equation*}
    dX^{\alpha}_t=b(t,X^{\alpha}_t,\nu^{\alpha}_t,\alpha_t)dt+\sigma(t,X^{\alpha}_t,\nu^{\alpha}_t,\alpha_t) dW_t,\quad
    X^{\alpha}_0=\xi.
    \end{equation*}
where $\nu^{\alpha}_t$ is a shorthand notation for $\mathcal{L}(X^{\alpha}_t,\alpha_t)$ and $\mu^{\alpha}_T$ is its first marginal at terminal time $T$. The dynamics of $X$ involves the law of this process, hence the terminology McKean-Vlasov dynamics~\citep{mckean1966class}. To alleviate notation we will sometimes write $\nu^* = \nu^{\alpha^*}$ for the law of the optimally controlled process.

\begin{remark}
    Although the two problems look similar, they in general have different solutions, \textit{i.e.}, $\hat\alpha \neq \alpha^*$ and $\hat\nu \neq \nu^*$, even when the functions in the cost and the dynamics are the same, see \textit{e.g.}~\citep{MR3968548}.
\end{remark}

\begin{remark}
    Although the mean field paradigm is the same, the special case where the interactions are only through the state distribution (\textit{i.e.}, the first marginal of $\nu$) has attracted more interest in the literature than the present general setup. However interactions through the distribution of controls appears in many applications, particularly in economics and finance as already mentioned in the introduction. See next sections for some examples.
\end{remark}

\begin{remark}
    Although the reinforcement learning literature typically focuses on infinite horizon discounted problems, we focus here on finite horizon problems. This will cause some numerical difficulties but is crucial for many applications. 
\end{remark}

\section{Two timescale approach}\label{sec:twotimescale}

\subsection{Discrete formulation}

To simplify the presentation and to be closer to the traditional reinforcement learning setup, we consider a discrete time model with a finite number of states and actions. Let $\mathcal{X}$ and $\mathcal{A}$ be finite sets corresponding to spaces of states and actions respectively, which can correspond to discretized version of the continuous spaces (possibly after a truncation) used in the previous section. We denote by $\Delta^{|\mathcal{X}|}$ the simplex in dimension $|\mathcal{X}|$, which we identify with the space of probability measures on $\mathcal{X}$.  $\Delta^{|\mathcal{X}|\times|\mathcal{A}| }$ is defined similarly on the product space $\mathcal{X}\times\mathcal{A}$. Moreover, we consider a discrete time setting, which here again could come from a suitable discretization of the continuous time evolution used in the previous section. We take a uniform grid, say $t_n = n \times \Delta t$, $n=0,1,2,\dots, N_T$, where $\Delta t = T/N_T>0$ is the time step. Hence $X_{t_n}$ and $\alpha_{t_n}$ in this section can be interpreted as approximation of the state and the action at time $t_n$ in the previous section. The state follows a random evolution in which $X_{t_{n+1}}$ is determined as a function of the current state $X_{t_n}$, the action $\alpha_{t_n}$, the state-action population distribution $\nu_{t_n}$ at time $t_n$, and some noise. We introduce the transition probability function:
$$
    p(x'|x,a,\nu), \qquad (x, x', a, \nu) \in \mathcal{X} \times \mathcal{X} \times \mathcal{A} \times \Delta^{|\mathcal{X} \times \mathcal{A}|}, 
$$
which gives the probability to jump to state $x'$ when being at state $x$ and using action $a$ and when the population distribution is $\nu$. For simplicity, we consider the homogeneous case where this function does not depend on time, which corresponds, in the continuous formulation, to the case where both $b$ and $\sigma$ are time-independent.
Restoring this time-dependence if needed is a straightforward procedure.

We now consider the MFG cost function given by: for $\nu = (\nu_{t_n})_{n=0,\dots,N_T}$ 
$$
    \widetilde J_{\nu}(\alpha) = \mathbb{E}\left[ \sum_{n=0}^{N_T-1} f(X^{\alpha,\nu}_{t_n}, \alpha_{t_n}, \nu_{t_n}) + g(X^{\alpha,\nu}_{t_{N_T}}, \mu_{t_{N_T}}) \right],
$$
where 
$\mu_{t_{N_T}}$ is the first marginal of $\nu_{t_{N_T}}$. Again, for simplicity, we assume that $f$ doesn't depend on time. The process $X^{\alpha,\nu}$ has a given initial distribution $\mu_0 \in \Delta^{|\mathcal{X}|}$ and follows the dynamics 
$$
    \mathbb{P}(X^{\alpha,\nu}_{t_{n+1}} = x' | X^{\alpha,\nu}_{t_n} = x, \alpha_{t_n} = a, \nu_{t_n} = \nu) = p(x'|x,a,\nu).
$$

Given a population distribution sequence $\nu$, the value function of an infinitesimal player is
$$
    V_\nu(x) = \inf_{\alpha}  V_\nu^\alpha(x),
$$
where
$$
    V_\nu^{\alpha}(x) = \mathbb{E}\left[ \sum_{n=0}^{N_T-1} f(X^{\alpha,\nu}_{t_n}, \alpha_{t_n}, \nu_{t_n}) + g(X^{\alpha,\nu}_{t_{N_T}}, \mu_{t_{N_T}})  \Big| X^{\alpha,\nu}_{0} = x\right].
$$
Note that the $\widetilde{J}$ and the $V$ functions are related by:
$$
    \widetilde J_{\nu}(\alpha) = \mathbb{E}_{X_0 \sim \mu_0} [V^{\alpha}_\nu(X_0)]. 
$$

On the other hand, we also consider the MFC cost function
$$
    \widetilde J(\alpha) = \mathbb{E} \left[ \sum_{n=0}^{N_T-1} f(X^{\alpha}_{n}, \alpha_{t_n}, \nu^\alpha_n) + g(X^{\alpha}_{t_{N_T}}, \mu^\alpha_{t_{N_T}}) \right],
$$
where $\nu^{\alpha}_{t_n} = \mathcal{L}(X^{\alpha}_{t_n},\alpha_{t_n})$ is the state-action distribution at time $t_n$ of $X^{\alpha}$ controlled by $\alpha$. The process $X^{\alpha}$ has initial distribution $\mu_0$ and dynamics 
$$
    \mathbb{P}(X^{\alpha}_{t_{n+1}} = x' | X^{\alpha}_{t_n} = x, \alpha_{t_n} = a, \nu_{t_n} = \nu^{\alpha}_{t_n}) = p(x'|x,a,\nu^{\alpha}_{t_n}).
$$ 
Since the dynamics is of MKV type, in general, the value function in a MFC problem is the value function of the social planner and it takes the distribution $\nu$ as input, see \textit{e.g.}~\citep{MR3258261_DPMFC_CRAS,pham2016discrete,CarmonaLauriereTan-2019-MQFRL,motte2019mean,gu2019dynamic,djete2019mckean}.  However, when the population is already evolving according to the sequence of distributions $\nu^{\alpha}$ generated by a control $\alpha$, the cost-to-go of an infinitesimal agent starting at position $x$ and using control $\alpha$ too is simply a function of its position and is given by
$$
    V^{\alpha}(x)   
    = \mathbb{E} \left[ \sum_{n=0}^{N_T-1} f(X^{\alpha}_{t_n}, \alpha_{t_n}, \nu^\alpha_{t_n}) + g(X^{\alpha}_{t_{N_T}}, \mu^\alpha_{t_{N_T}}) \Big| X^{\alpha}_0 = x\right].
$$

\subsection{Action-value function}

The state value function is useful as far as the value of the game or control problem is concerned. However, it does not provide any information about the equilibrium or optimal control $\hat\alpha$ or $\alpha^*$. For this reason, one can introduce the state-action value function, also called $Q$-function, which takes as inputs not only a state $x$ but also an action $a$. For a standard Markov Decision Process (MDP without mean field interactions) the $Q$-function characterizes the optimal cost-to-go if one starts at state $x$ and uses action $a$ before starting using the optimal control. To approximate this function, one of the most popular methods in RL is the so-called Q-learning~\citep{watkins1989learning}. See \textit{e.g.}~\cite[Chapter 3]{sutton2018reinforcement} for more details.

Before moving on to the mean-field setup, let us recall that in the traditional setup, the definition of the optimal $Q$-function, denoted by $Q^*$, is given by:
 \begin{equation*}%
 \left\{
 \begin{aligned}
 &Q^*_{N_T}(x,a) = g(x), \qquad (x,a) \in \mathcal{X} \times \mathcal{A},
 \\
    &Q^*_n(x,a)=\min_\alpha\EE\left[\sum_{n'=n}^{N_T-1}  f(X_{t_{n'}},\alpha_{n'}(X_{t_{n'}})) + g(X_{t_{N_T}})\,\Big\vert\, X_{t_n}=x,A_{t_n}=a\right], 
    \\
    &\qquad\qquad n < N_T, (x,a) \in \mathcal{X} \times \mathcal{A} .
    \end{aligned}
    \right.
\end{equation*}
where $\alpha_{n'}(\cdot) = \alpha(t_{n'},\cdot)$. Using dynamic programming, it can be shown that $(Q^*_n)_n$ is the solution of the Bellman equation:  
 \begin{equation*}%
 \left\{
 \begin{aligned}
 &Q^*_{N_T}(x,a) = g(x), \qquad (x,a) \in \mathcal{X} \times \mathcal{A},
    \\
   & Q^*_n(x,a) = f(x, a) + \sum_{x' \in \mathcal{X}} p(x' | x, a) \min_{a'} Q^*_{n+1}(x',a'), \qquad n<N_T, (x,a) \in \mathcal{X} \times \mathcal{A}.
    \end{aligned}
    \right.
\end{equation*}

The corresponding optimal value function $(V_n^*)_n$ is given by:
$$
    V^*_n(x) = \min_a Q^*_n(x,a), \qquad n \le N_T, x \in \mathcal{X}.
$$
As mentioned above, one of the main advantages of computing the action-value function instead of the value function is that from the former, one can directly recover the optimal control at time $n$, given by $\argmin_{a \in \mathcal{A}} Q^*_n(x,a)$. This is particularly important in order to design model-free methods, as we will see in the next section.

The above approach can be adapted to solve MFG by noticing that, when the population behavior is given, the problem posed to a single representative agent is a standard MDP. It can thus be tackled using a $Q$-function which implicitly depends on the population distribution: given $\nu = (\nu_{t_n})_{n = 0,\dots,N_T}$
 \begin{equation*}%
 \left\{
 \begin{aligned}
 &Q^*_{N_T,\nu}(x,a) = g(x,\mu_{t_{N_T}}), \qquad (x,a) \in \mathcal{X} \times \mathcal{A},
    \\
   & Q^*_{n,\nu}(x,a) = f(x, a, \nu_{t_n}) 
   \\
   &\qquad + \sum_{x' \in \mathcal{X}} p(x' | x, a,\nu_{t_n}) \min_{a'} Q^*_{n+1,\nu}(x',a'), \qquad n<N_T, (x,a) \in \mathcal{X} \times \mathcal{A}.
    \end{aligned}
    \right.
\end{equation*}
This function characterizes, at each time step $t_n$, the optimal cost-to-go for an agent starting at time $t_n$ at state $x$, using action $a$ for the first step, and then acting optimally for the rest of the time steps, while the population evolution is given by $\nu = (\nu_{t_n})_n$. However, to find the Nash equilibrium, it is not sufficient to compute the $Q$-function for an arbitrary sequence of distributions $\nu$: we want to find $Q^*_{\nu^*}$ where $\nu^*$ is the population evolution generated by the optimal control computed from $Q^*_{\nu^*}$. In the sequel, we will directly aim at the $Q$-function $Q^*_{\nu^*}$ via a two timescale approach.

In the MFC problem the population distribution is not fixed while each player optimizes because all the agents cooperate to choose a distribution which is optimal from the point of view of the whole society. As a consequence, the optimization problem can not be recast as a standard MDP. However we will show below that it is still possible to compute the social optimum using a modified $Q$-function (not involving explicitly the population distribution). This major difficulty is treated in detail in the context of infinite horizon in our previous work \citep{angiuli2020unified}.

\subsection{Unification through a two timescale approach}\label{sec: unified 2scale approach}

A simple approach to compute the MFG solution is to iteratively update the state-action value function, $Q$, and the population distribution, $\nu$: Starting with an initial guess $\nu^{(0)}$, repeat for $k=0,1,\dots$,
\begin{enumerate}
    \item Solve the backward equation for $Q^{(k+1)} = Q^{*}_{\nu^{(k)}}$, which characterizes the optimal state-action value function of a typical player if the population behavior is given by $\nu^{(k)}$:
     \begin{equation}%
     \label{eq:backward-Q-kp1}
         \left\{
         \begin{aligned}
         &Q^{(k+1)}_{N_T}(x,a) = g(x,\mu^{(k)}_{t_{N_T}}), \qquad (x,a) \in \mathcal{X} \times \mathcal{A},
            \\
           & Q^{(k+1)}_{n}(x,a) = f(x, a,\nu^{(k)}_{t_n}) 
           \\
           &\qquad + \sum_{x' \in \mathcal{X}} p(x' | x, a,\nu^{(k)}_{t_n}) \min_{a'} Q^{(k+1)}_{n+1}(x',a'), \quad  n<N_T, (x,a) \in \mathcal{X} \times \mathcal{A}.
            \end{aligned}
            \right.
    \end{equation}
    \item Solve the forward equation for $\mu^{(k+1)}$ (resp. $\nu^{(k+1)}$), which characterizes the evolution of the population state distribution (resp. state-action distribution) if everyone uses the optimal control $\alpha^{(k+1)}_{t_n}(x) = \argmax_a Q^{(k+1)}_{n}(x,a)$ coming from the above $Q$-function  (assuming this control is uniquely defined for simplicity):
        \begin{equation}%
        \label{eq:forward-mu-nu-kp1}
         \left\{
         \begin{aligned}
         &\mu^{(k+1)}_{t_0}(x) = \mu_{t_0}(x), \qquad x \in \mathcal{X},
         \\
         &\nu^{(k+1)}_{t_0}(x,a) = \mu_{t_0}(x) \mathbf{1}_{a = \alpha^{(k+1)}_{t_n}(x)}, \qquad (x,a) \in \mathcal{X} \times \mathcal{A},
            \\
           & \mu^{(k+1)}_{t_{n+1}}(x) = \sum_{x' \in \mathcal{X}} p(x' | x, \alpha^{(k+1)}_{t_n}(x'),\nu^{(k+1)}_{t_n}) , \qquad 0 \le n < N_T, x \in \mathcal{X},
            \\
           & \nu^{(k+1)}_{t_{n+1}}(x,a) = \mu^{(k+1)}_{t_{n+1}}(x) \mathbf{1}_{a = \alpha^{(k+1)}_{t_{n+1}}(x)} , \qquad 0 \le n < N_T, (x,a) \in \mathcal{X} \times \mathcal{A}.
            \end{aligned}
            \right.
    \end{equation}
\end{enumerate}
Here the evolution of the joint state-action population distribution is simply the product of the state distribution and a Dirac mass:
$$
    \nu^{(k+1)}_{t_n} = \mu^{(k+1)}_{t_n} \otimes \delta_{\alpha^{(k+1)}_{t_n}}.
$$
This is because we assumed that the optimal control is given by a deterministic function from $\mathcal{X}$ to $\mathcal{A}$. If we were using randomized control, the Dirac mass would need to be replaced by the distribution of controls. 

To alleviate notation, let us introduce the operators $\widetilde{\mathcal{T}}: (\Delta^{|\mathcal{X}\times\mathcal{A}|})^{N_T+1} \to (\mathbb{R}^{|\mathcal{X}\times\mathcal{A}|})^{N_T+1}$ and $\widetilde{\mathcal{P}}: (\mathbb{R}^{|\mathcal{X}\times\mathcal{A}|})^{N_T+1} \to (\Delta^{|\mathcal{X}\times\mathcal{A}|})^{N_T+1}$  such that: \eqref{eq:backward-Q-kp1} and~\eqref{eq:forward-mu-nu-kp1} rewrite
$$
    Q^{(k+1)} = \widetilde{\mathcal{T}}(\nu^{(k)}),
    \qquad
    \nu^{(k+1)} = \widetilde{\mathcal{P}}(Q^{(k+1)}).
$$
If this iteration procedure converges, we have $Q^{(k+1)} \to Q^{(\infty)}, \nu^{(k+1)} \to \nu^{(\infty)}$ as $k \to +\infty$ for some $Q^{(\infty)}, \nu^{(\infty)}$ satisfying
$$
    Q^{(\infty)} = \widetilde{\mathcal{T}}(\nu^{(\infty)}),
    \qquad
    \nu^{(\infty)} = \widetilde{\mathcal{P}}(Q^{(\infty)}),
$$
which implies that $\nu^{(\infty)}$ is the state-action equilibrium distribution of the MFG solution, and the associated equilibrium control is given by: $\alpha^{(\infty)}_{t_n}(x) = \argmax_a Q^{(\infty)}_{n}(x,a)$ for each $n$. 

However, this procedure fails to converge in many MFG by lack of strict contraction property. To remedy this issue, a simple twist is to introduce some kind of damping. Building on this idea, we introduce the following iterative procedure, where $(\rho_Q^{(k)})_{k \ge 0}$ and $(\rho_\nu^{(k)})_{k \ge 0}$ are two sequences of learning rates:
$$
    Q^{(k+1)} = (1-\rho_Q^{(k)})Q^{(k)} + \rho_Q^{(k)} \widetilde{\mathcal{T}}(\nu^{(k)}),
    \qquad
    \nu^{(k+1)} = (1-\rho_\nu^{(k)})\nu^{(k)} + \rho_\nu^{(k)} \widetilde{\mathcal{P}}(Q^{(k+1)}).
$$
For the sake of brevity, let us introduce the operators $\mathcal{T}: (\mathbb{R}^{|\mathcal{X}\times\mathcal{A}|})^{N_T+1} \times (\Delta^{|\mathcal{X}\times\mathcal{A}|})^{N_T+1} \to (\mathbb{R}^{|\mathcal{X}\times\mathcal{A}|})^{N_T+1}$ and $\mathcal{P}: (\mathbb{R}^{|\mathcal{X}\times\mathcal{A}|})^{N_T+1} \times (\Delta^{|\mathcal{X}\times\mathcal{A}|})^{N_T+1} \to (\Delta^{|\mathcal{X}\times\mathcal{A}|})^{N_T+1}$ 
$$
    \mathcal{T}(Q, \nu) = \widetilde{\mathcal{T}}(\nu) - Q,
    \qquad
    \mathcal{P}(Q, \nu) = \widetilde{\mathcal{P}}(Q) - \nu.
$$
Then the above iterations can be written as
\begin{equation}
\label{eq:2scale-Q-nu-k}
    Q^{(k+1)} = Q^{(k)} + \rho_Q^{(k)} \mathcal{T}(Q^{(k)}, \nu^{(k)}),
    \qquad
    \nu^{(k+1)} = \nu^{(k)} + \rho_\nu^{(k)} \mathcal{P}(Q^{(k+1)}, \nu^{(k)}).
\end{equation}

If $\rho^{(k)}_\nu < \rho^{(k)}_Q$, the $Q$-function is updated at a faster rate, while it is the converse if $\rho^{(k)}_\nu > \rho^{(k)}_Q$. We can thus intuitively guess that these two regimes should converge to different limits. Similar ideas have been studied by Borkar~\citep{borkar1997stochastic,MR2442439} in the so-called two timescales approach. The key insight comes from rewriting the (discrete time) iterations in continuous time as a pair of ODEs. From~\cite[Chapter 6, Theorem 2]{MR2442439}, we expect to have the following two situations:

\begin{itemize}
    \item If $\rho^{(k)}_\nu < \rho^{(k)}_Q$, the system \eqref{eq:2scale-Q-nu-k} tracks the ODE system
    \begin{subequations}%
     \begin{empheq}[left=\empheqlbrace]{align*}
    \dot Q^{(t)} &= \frac{1}{\epsilon} \mathcal{T}(Q^{(t)}, \nu^{(t)}),
    \\
    \dot \nu^{(t)} &= \mathcal{P}(Q^{(t)}, \nu^{(t)}),
    \end{empheq}
    \end{subequations}
    where $\rho^{(k)}_\nu /\rho^{(k)}_Q$  is thought of being of order $\epsilon\ll 1$.
    Hence, for any fixed $\tilde \nu$, the solution of 
    $$
    \dot Q^{(t)} = \frac{1}{\epsilon} \mathcal{T}(Q^{(t)}, \tilde \nu),
    $$ 
    is expected to converge as $\epsilon \to 0$ to a $Q^{\tilde \nu}$ such that $\mathcal{T}(Q^{\tilde \nu}, \tilde \nu) = 0$. This condition can be interpreted as the fact that $Q^{\tilde \nu} = (Q^{\tilde \nu}_{t_n})_{n=0,\dots,N_T}$ is the state-action value function of an infinitesimal agent facing the crowd distribution sequence $\tilde \nu = (\tilde \nu_{t_n})_{n=0,\dots,N_T}$. Then the second ODE becomes 
    $$
        \dot \nu^{(t)} = \mathcal{P}(Q^{\nu^{(t)}}, \nu^{(t)}),
    $$
    which is expected to converge as $t \to +\infty$ to a $\nu^{(\infty)}$ satisfying
    $$
        \mathcal{P}(Q^{\nu^{(\infty)}}, \nu^{(\infty)}) = 0.
    $$
    This condition means that $\nu^{(\infty)}$ and the associated control given by $\hat\alpha_n(x) = \argmin_{a} Q^{\nu^{(\infty)}}_{n}(x, a)$ form a Nash equilibrium.
    
    \item If $\rho^{(k)}_\nu > \rho^{(k)}_Q$, the system \eqref{eq:2scale-Q-nu-k} tracks the ODE system
    \begin{subequations}%
     \begin{empheq}[left=\empheqlbrace]{align*}
    \dot Q^{(t)} &= \mathcal{T}(Q^{(t)}, \nu^{(t)}),
    \\
    \dot \nu^{(t)} &= \frac{1}{\epsilon}  \mathcal{P}(Q^{(t)}, \nu^{(t)}),
    \end{empheq}
    \end{subequations}
     where $\rho^{(k)}_Q /\rho^{(k)}_\nu$  is thought of being of order $\epsilon\ll 1$.
    Here, for any fixed $\tilde Q$, the solution of 
    $$
        \dot \nu^{(t)} = \frac{1}{\epsilon}  \mathcal{P}(\tilde Q, \nu^{(t)}),
    $$
     is expected to converge as $\epsilon \to 0$ to a $\nu^{\tilde Q}$ such that $\mathcal{T}(\tilde Q, \nu^{\tilde Q}) = 0$, meaning that $\nu^{\tilde Q} = (\nu^{\tilde Q}_{t_n})_{n=0,\dots,N_T}$ is the  distribution evolution of a population in which every agent uses control $\tilde \alpha_{t_n}(x) = \argmin_{a} \tilde Q_n(x,a)$ at time $t_n$. In fact, the definitions of 
     $\tilde \alpha_{t_n}$ and $\nu^{\tilde Q}$ need to be {\it modified} to take into account the first action $(x,a)$. For the details of this crucial step for handling MFC, we refer to \citep{angiuli2020unified}.
     
    Then the first ODE becomes 
    $$
       \dot Q^{(t)} = \frac{1}{\epsilon} \mathcal{T}(Q^{(t)}, \nu^{Q^{(t)}}),
    $$
    which is expected to converge as $t \to +\infty$ to a  $Q^{(\infty)}$ such that
    $$
        \mathcal{T}(Q^{(\infty)}, \nu^{Q^{(\infty)}}) = 0.
    $$
    This condition (in the {\it modified} MFC setup) means that the control $\hat\alpha_{t_n}(x) = \argmin_{a} Q^{(\infty)}_{n}(x, a)$ is a MFC optimum and the induced optimal distribution is $\nu^{Q^{(\infty)}}$.
\end{itemize}

The above iterative procedure is purely deterministic and allows us to understand the rationale behind the two timescale approach. However, in practice we rarely have access to the operators $\mathcal{T}$ and $\mathcal{P}$. Instead, we will consider that we only have access to noisy versions and we use intuition from stochastic approximation to design an algorithm. Instead of assuming that we know the dynamics or the reward functions, we will simply assume that the learning agent can interact with an environment from which she can sample stochastic transitions.

\section{Reinforcement learning algorithm}\label{sec:RLalgo}

\subsection{Reinforcement learning}

RL is a branch of Machine Learning which studies algorithms to solve a MDP based on trials and errors. An MDP describes the sequential interaction of an agent with an environment. Let $\mc{X}$ and $\mc{A}$ be the state and action space respectively. At each time $t_n$, the agent observes its current state $X_{t_n} \in \mathcal{X}$ and chooses an action $A_{t_n} \in \mathcal{A}$. Due to the agent's action, the environment provides the new state of the agent $X_{t_{n+1}}$ and incurs a cost $f_{t_{n+1}}$. The goal of the agent is to find an optimal strategy (or policy) $\pi^*$ which assigns to each state an action in order to minimize the aggregated discounted costs. The aim of RL is to design methods which allow the agent to learn (an approximation of) $\pi^*$ by making repeated use of the environment's outputs but without knowing how the environment produces the new state and the associated cost. A detailed overview of this field can be found in \cite{sutton2018reinforcement} (although RL methods are often presented with reward maximization objectives, we consider cost minimization problems for the sake of consistency with the MFG literature). 

Here and in what follows, we use {\it policy} $\pi$ instead of control $\alpha$ as the algorithm uses in fact $\epsilon$-greedy policies. In the limit $\epsilon\to 0$ the optimal policy is in fact a deterministic control.

\subsection{Algorithm}\label{sec:algorithm}

In this section we propose an extension of the Unified Two Timescales Mean Field Q-learning (U2-MF-QL) algorithm discussed in our previous work \citep{angiuli2020unified}. 

Q-learning is one of the most popular procedure in RL introduced by Watkins in his seminal work \citep{watkins1989learning}. It is designed to solve problems with finite and discrete state and action spaces, $\mc{X}$ and $\mc{A}$. It is based on the evaluation of the optimal action-value function $Q^*(x,a)$, defined in the case of an infinite horizon minimization problem as
 \begin{align*}%
    Q^*(x,a)=\min_\pi\EE\left[\sum_{n=0}^\infty \rho^n f_{t_{n+1}}(X_{t_n},\pi(X_{t_n})) \,\Big\vert\, X_{t_0}=x,A_{t_0}=a\right],
\end{align*}
 which represents the optimal expected aggregated discounted cost when starting in the state $x$ and choosing the first action to be $a$. The optimal action at state $x$ is provided by the argmin of $Q^*(x,\cdot)$, \textit{i.e.}, $\pi^*(x) = \argmin_{\mathcal{A}} Q^*(x,\cdot)$. However $Q^*$ is a priori unknown. In order to learn $Q^*$ by trials and errors, an approximate version $Q$ of the table $Q^*$ is constructed through a stochastic approximation procedure based on the Bellman equation given by 
\begin{equation}
Q^*(x,a)=     \EE\left[   f_{t_{0}}(X_{t_0},A_{t_0}) + Q^*(X_{t_1},\pi^*(X_{t_1})) \Big\vert\, X_{t_0}=x,A_{t_0}=a\right].
\end{equation}  At each step, an action is taken, which leads to a cost and to a new state. On the one hand, it is interesting to act efficiently in order to avoid high costs, and on the other hand it is important to improve the quality of the table $Q$ by trying actions and states which have not been visited many times so far. This is the so-called exploitation--exploration trade-off. The trade-off between exploration of the unknown environment and exploitation of the currently available information is taken care of by an $\epsilon$-greedy policy based on $Q$. The algorithm chooses the action that minimizes the immediate cost with probability $1-\epsilon$, and a random action otherwise. 

The U2-MF-QL algorithm represents a unified approach to solve asymptotic Mean Field Games and Mean Field Control problems based on the relationship between two learning rates relative to the update rules of the $Q$ table and the distribution of the population $\mu$ respectively. Based on the intuition presented in Section \ref{sec:twotimescale}, a choice of learning rates $(\rho^Q, \rho^{\mu})$ such that $\rho^Q > \rho^{\mu}$ allows the algorithm to solve a MFG problem. The estimation of $Q$ is updated at a faster pace with respect to the distribution which behaves as quasi-static mimicking the freezing of the flow of measures characteristic of the solving scheme discussed in Section \ref{sec:mfg}. On the other hand, learning rates satisfying $\rho^Q < \rho^{\mu}$ allow the algorithm to updates instantaneously the control function (Q table) at any change of the distribution reproducing the  MFC framework. Under suitable assumptions, one may expect the asymptotic problems to be characterized by controls that are independent of time. In this case, the learning goals reduce to a control function valid for every time point and the asymptotic distribution of the states of the population. \\
The finite horizon framework presented in Sections \ref{sec:mfg} and \ref{sec:mfc} differs from the asymptotic case discussed in \citep{angiuli2020unified} in several ways other than the restriction on the finite time interval $[0,T]$. First, the mean field interaction is through the joint distribution of states and actions of the population rather than the marginal distribution of the states. Further, both the control rule and the mean field distribution are generally time dependent. Due to these differences, the $2-$dimensional matrix $Q$ in U2-MF-QL is replaced by a $3-$dimensional matrix $\boldsymbol{Q} \coloneqq \left(Q_n(\cdot,\cdot) \right)_{n=0,\dots,N_T} = \left(Q(\cdot,\cdot,t_n) \right)_{n=0,\dots,N_T}$ in the finite horizon version of the algorithm (U2-MF-QL-FH). The extra dimension is introduced to learn a time dependent control function. \\
The Unified Two Timescales Mean Field Q-learning for Finite Horizon problems (U2-MF-QL-FH) is designed to solve problems with finite state and action spaces in finite and discrete time. 
\begin{algorithm}
\caption{Unified Two Timescales Mean Field Q-learning - Finite Horizon \label{algo:U2MFQL}} 
\begin{algorithmic}[1]
\REQUIRE $\tau = \{t_0=0, \dots, t_{N_T}=T \}$ with $t_0=0 < \dots < t_{N_T}=T$  : time steps, 
\\ $\mathcal{{X}} =\{ x_0, \dots, x_{|\mathcal{{X}}|-1}\}$ : finite state space, \\ 
$\mathcal{{A}} =\{ a_0, \dots, a_{|\mathcal{{A}}|-1}\}$ : finite action space, \\
$\mu_0$ : initial distribution of the representative player,\\
$\epsilon$ : factor related to the $\epsilon-$greedy policy,\\
$tol_{\nu}$, $tol_Q$ : break rule tolerances.
\STATE \textbf{Initialization}: $Q_{n}(\cdot,\cdot) \coloneqq Q(\cdot,\cdot,t_n) =0$ for all $(x,a)\in \mathcal{{X}}\times \mathcal{{A}}$, for all $ t_n \in \tau$,\\ $\nu^0_{t_n}=\frac{1}{|\mathcal{{X}}\times\mathcal{{A}}|} J_{|\mathcal{{X}}|\times |\mathcal{{A}}|}$ for $n=0, \dots, N_T$ where $J_{d\times m}$ is an $d\times m$ unit matrix

\FOR{each episode $k=1,2,\dots$}
\STATE \textbf{Observe} $X^k_{t_0} \sim \mu_0$ \\
\FOR{$ n \gets 0$ to $N_T$} 

\STATE \textbf{Choose action $A^k_{t_n}$} using the $\epsilon$-greedy policy derived from $Q_n^{k-1}(X^k_{t_n},\cdot)$  
\\
\STATE
{\color{blue}\textbf{Update} $\nu$: \\
$\nu^{k}_{t_{n}} = \nu^{k-1}_{t_{n}} + \rho^{\nu}_{k} (\boldsymbol{\delta}(X^k_{t_{n}},A^k_{t_{n}}) - \nu^{k-1}_{t_{n}})$\\ where $\boldsymbol{\delta}(X^k_{t_{n}},A^k_{t_{n}}) = \left(\mathbf{1}_{x,a}(X^k_{t_{n}},A^k_{t_{n}})\right)_{x\in \mc{X},a\in \mc{A}} $}
\\
\textbf{Observe cost $f_{t_{n+1}}=f(X^k_{t_n},A^k_{t_n},\nu^{k}_{t_n})$} and state $X^k_{t_{n+1}}$ provided by the environment %

\STATE
\color{red}\textbf{Update} $Q_n$:\\ $Q_n^k(x,a) \coloneqq $\\
$\coloneqq\begin{cases} 
Q_n^{k-1}(x,a)+\rho^{Q_n}_{x,a,k} [\mc{B} -Q_n^{k-1}(x,a)] &\text{ if } (X^k_{t_n},A^k_{t_n})=(x,a) \\
Q_n^{k-1}(x,a)  &\text{ o.w. } 
\end{cases}$\\
where \\
$\mc{B} \coloneqq \begin{cases} f_{t_{n+1}} +\gamma \min_{a'\in \mc{A}} Q_{n+1}^{k-1}(X^k_{t_{n+1}},a'), \quad &\text{if } t_n < T \\f_{t_{n+1} } , \quad &\text{o.w.} \end{cases} $\\
\ENDFOR
\IF{ $ \lvert|{\nu_{t_n}^k-\nu_{t_n}^{k-1}}\rvert|_1 \leq tol_{\nu}$ and $\|Q_n^k-Q_n^{k-1}\|_{1,1}<tol_Q$ for all $n=0,\dots,N_T$}
\STATE break
\ENDIF
\ENDFOR
\end{algorithmic}
\end{algorithm}
The same algorithm  can be applied to MFG and MFC problems where the the interaction with the population is through the marginal distribution of the states $\mu \in \mathcal{P}(\mc{X})$ or the law of the controls $\theta \in \mathcal{P}(\mc{A})$. In these cases the estimation of the flow of marginal distributions is obtained through the vectors $\left( \mu_{t_n} \right)_{n=0,\dots,N_T}$ (resp. $\left( \theta_{t_n} \right)_{n=0,\dots,N_T}$)  defined on the space $\mc{X}$ (resp. $\mc{A}$). The initialization is given by $\mu^0_{t_n}=\left[ \frac{1}{|\mathcal{{X}}|}, \dots, \frac{1}{|\mathcal{{X}}|} \right] $ $\left(\text{resp. } \theta^0_{t_n}=\left[ \frac{1}{|\mathcal{{A}}|}, \dots, \frac{1}{|\mathcal{{A}}|} \right] \right) $ for $n=0, \dots, N_T$. The update rule at episode $k$ is given by $\mu^{k}_{t_{n}} = \mu^{k-1}_{t_{n}} + \rho^{\mu}_{k} (\boldsymbol{\delta}(X_{t_{n}}) - \mu^{k-1}_{t_{n}})$  $\left( \text{resp. }\theta^{k}_{t_{n}} = \theta^{k-1}_{t_{n}} + \rho^{\theta}_{k} (\boldsymbol{\delta}(A_{t_{n}}) - \theta^{k-1}_{t_{n}}\right)$ where $\boldsymbol{\delta}(X_{t_{n}}) = \bracket{\mathbf{1}_{x_0}(X_{t_{n}}), \dots, \mathbf{1}_{x_{|{\mc{X}}|-1}}(X_{t_{n}})}$  $\left(\text{resp. }\boldsymbol{\delta}(A_{t_{n}}) = \bracket{\mathbf{1}_{a_0}(A_{t_{n}}), \dots, \mathbf{1}_{a_{|{\mc{A}}|-1}}(A_{t_{n}})}\right)$ for $n=0, \dots, N_T$.\\
\subsection{Learning rates}
The algorithm \ref{algo:U2MFQL} is based on two stochastic approximation rules for the distribution $\boldsymbol{\nu}$ and the $3-$dim matrix $\boldsymbol{Q}$. The design of the learning is discussed widely in the literature, in a general context by  \citep{borkar1997stochastic} and \citep{MR2442439}, and with focus in reinforcement learning by \citep{borkar1997actor} and  \citep{even2003learning}.
Based on experimental evidences, we define the learning rates appearing in Algorithm~\ref{algo:U2MFQL} as follows: 
\begin{equation} \label{eqn : rates}
 \rho^{Q_{n}}_{x,a,k}=\frac{1}{\left(1+N_T\#|(x,a,t_n,k)|\right)^{\omega^Q}},     \quad \quad \rho^{\nu}_k=\frac{1}{(1+k)^{\omega^{\nu}}}, 
\end{equation}
 where $\#|(x,a,t_n,k)|$ is counting the number of visits of the pair $(x,a)$ at time $t_n$ until episode $k$. Differently from the asymptotic version of the algorithm presented in \citep{angiuli2020unified} for which each pair $(x,a)$ has a unique counter for all time points, in the finite horizon formulation a distinct counter $\#|(x,a,t_n,k)|$ is defined for each time point $t_n$. This choice of learning rates allows to update each matrix $Q_{n}$  in an asynchronous way. The exponent $\omega^Q$ can take values in $(\frac{1}{2},1].$ As presented in Section~\ref{sec: unified 2scale approach}, the pair $(\omega^Q,\omega^{\nu})$ is chosen depending on the particular problem to solve. In a competitive framework (MFG), these parameters have to be searched in the set of values for which the condition $\rho^Q>\rho^{\nu}$ is satisfied at each iteration. On the other hand, a good choice for the cooperative case (MFC) should satisfy the condition $\rho^Q<\rho^{\nu}$.
 
 \subsection{Application to continuous problems}\label{sec:continuous_extension_algo}

Although it is presented in a setting with finite state and action spaces, the application of the algorithm U2-MF-QL-FH can be extended to continuous problems. Such adaptation requires truncation and discretization procedures to time, state and action spaces which should be calibrated based on the specific problem.

In practice, a continuous time interval $[0,T]$ would be replaced by a uniform discretization $\tau = \{t_n \}_{n\in \{0,\dots, N_T\}} $. The environment would provide the new state and reward at these discrete times. The continuous state would be projected on a finite set $\mc{X} = \{ x_0, \dots, x_{|\mc{X}|-1} \} \subset \RR^d$. Likewise, actions will be provided to the environment in a finite set $ \mc{A} = \{ a_0, \dots, a_{|\mc{A}|-1} \} \subset \RR^k $, where the projected distribution $\nu$  would be estimated. Then Algorithm~\ref{algo:U2MFQL} is ran on those spaces.

In the problems presented in Section~\ref{sec:trader}, we will use the benchmark linear-quadratic models given in continuous time and space for which we present explicit formulas.
In that case, we use an Euler discretization of the dynamics followed by a projection on $\mc{X}$. We do not address here the error of approximation since the purpose of this comparison with a benchmark is mainly for illustration.

\section{A mean field accumulation problem}\label{sec:accumulation}

\subsection{Description of the problem}

A further application of mean field theory to economics is given by the mean field capital accumulation problem by Huang in \citep{huang2013mean}. In this paper, the author studies an extension of the the classical one-agent modeling of optimal stochastic growth to an infinite population of symmetric agents. We introduce the model following the author's presentation.   \newline

\noindent At discrete time $t \in \mathbb{Z}_+$, the wealth of the representative agent is represented by a process $X_t^{\balpha,\btheta}$ characterized by the dynamics
\begin{equation}
X_{t+1}^{\balpha,\btheta} = G\left(\int a d{\theta}_t(a), W_t\right) \alpha_t
\end{equation}
where ${\balpha = (\alpha_t)_{0 \leq t \leq T}}$ is the  controlled variable denoting the agent's investment for production, $  G\left(\int a d{\theta}_t(a), W_t\right)$ is the production function,  $\btheta = (\theta_t)_{0 \leq t \leq T}$ is the mean field term represented by the law of the investment level of the population, $  \int a d{\theta}_t(a)$ is its mean, and  $\bg{W} = (W_t)_{0 \leq t \leq T}$ is a random disturbance. At each time $t$, the control ${\alpha_t}$ can only take values in $[0,X_t^{\balpha,\btheta}]$ so that $Supp(\theta_t) \subseteq [0,X_t^{\balpha,\btheta}]$,  implying that borrowing is not allowed. The wealth remaining after investment is all consumed, \textit{i.e.} the consumption variable $c_t$ is equal to $ c_t = X_t^{\balpha,\btheta} - \alpha_t$. The model is based on the following assumptions:
\begin{itemize}
    \item [(A1)] $\bg{W}$ is a random noise source with support $D_W$. The initial state $X_{t_0}$ is a positive random variable independent of $\bg{W}$ with mean $m_0$; 
    \item [(A2)] The function $G : [0,\infty) \times D_W \mapsto  [0,\infty)$ is continuous. If $w \in D_W$ is fixed, $G(z,w)$ is a decreasing function of $z$;
    \item[(A3)] $\mathbb{E}G(0,W)< \infty$ and $\mathbb{E}G(z,W)> 0 $ for each $z \in [0,\infty)$. 
\end{itemize}
The multiplicative factor $G$ in the dynamics of the wealth process $X_t^{\balpha,\btheta}$ shows the direct dependence of the wealth on  both the individual investment and the population aggregated investment. Further, assumption (A2) relates to the negative mean field impact explained as the loss in production efficiency when the aggregated investment increases. An example for the function $G$ is given by $G(z,w) = \frac{\beta w }{ 1 + \delta z^{\eta}},$ where $\beta $, $\delta $, $\eta$ are non negative parameters. Let $\bg{W}$ be a positive random noise with mean equal to 1. Then $D_W \subset [0, \infty)$ and (A2) - (A3) are satisfied. \newline

\noindent The goal of the agent is to optimize the expected aggregated discounted utility of consumption given by
\begin{equation}
    J(\balpha,\btheta) = \mathbb{E} \sum_{t=0}^T \rho^t v(c_t) = \mathbb{E} \sum_{t=0}^T \rho^t v(X_t^{\balpha,\btheta}-\alpha_t),
\end{equation}
where $\rho \in (0,1]$ is the discount factor. In particular, the author of \citep{huang2013mean} analyses the case of a Hyperbolic Absolute Risk Aversion (HARA) utility function defined as
\begin{equation}
 \label{eq:hara_utility}
v(c_t) = v(X_t^{\balpha,\btheta}-\alpha_t) \coloneqq \frac{1}{\gamma} (X_t^{\balpha,\btheta}-\alpha_t)^{\gamma},   
\end{equation}
where $\gamma \in (0,1)$. \newline

\subsection{Solution of the MFG}\label{sec:hara_mfg_sol}

\noindent In a competitive game setting, the resulting mean field game problem has solution given by Theorem~3 of Section~3.2 and Theorem~6 of Section~4 in \citep{huang2013mean}. Let denote the functions $\Phi(z)$, $\phi(z)$ and $\Psi(z)$ as follows
\[ \Phi(z) = \rho \mathbb{E} G^{\gamma} (z,W), \quad  \phi(z) = \Phi(z)^{\frac{1}{\gamma - 1 }}, \quad  \Psi(z) = \mathbb{E} G(z,W).\]
Let suppose that the mean field interaction is through $(z_t)_{t=0, \dots, T}$ the first moment of the flow of measures $\btheta = (\theta_t)_{t=0, \dots, T}$. The relative value function is defined as
\[ V^{\btheta}(t,x) = \sup_{{\balpha}} \mathbb{E} \left[ \sum_{s=t}^T \rho^{s} v(X_{s}^{\balpha,\btheta}-\alpha_{s}) | X_t^{\balpha,\btheta} = x \right].\] 
The value function is equal to  $V^{\btheta}(t,x) = \frac{1}{\gamma} D_t^{\gamma-1} x^{\gamma}$, where $D_t$ can be obtained using the recursive formula
\[ D_t = \frac{\phi(z_t) D_{t+1}}{1+\phi(z_t) D_{t+1}}, \quad  \quad D_T =1. \] The optimal control w.r.t. $\bg{\theta}$ is given by 
\[ \hat{\alpha}_t (x)= \frac{x}{1+\phi(z_t) D_{t+1}}, \quad t \leq T-1, \quad \quad \hat{\alpha}_T = 0.\]
The equivalent of the Nash equilibrium in the mean field limit is obtained by solving the fixed point equation 
\[ (\Lambda_0,\dots,\Lambda_{T-1})( z_0, \dots, z_{T-1}) = ( z_0, \dots, z_{T-1}), \]
 where

 \begin{equation*}
     \begin{cases} 
     \Lambda_0 ( z_0, \dots, z_{T-1}) \coloneqq \frac{1+ \phi(z_{T-1}) + \dots +\phi(z_{T-1}) \dots \phi(z_{1}) }{1+ \phi(z_{T-1}) + \dots +\phi(z_{T-1}) \dots \phi(z_{0})}m_0, \\
     \quad \\
     \Lambda_k ( z_0, \dots, z_{T-1}) \coloneqq\\
     \coloneqq \frac{1+ \phi(z_{T-1}) + \dots +\phi(z_{T-1}) \dots \phi(z_{k+1}) }{1+ \phi(z_{T-1}) + \dots +\phi(z_{T-1}) \dots \phi(z_{0})} \Psi(z_{k-1}) \dots \Psi(z_0) m_0, \quad \text{ for } 1 \leq k \leq T - 2, \\
      \quad \\
     \Lambda_{T-1} ( z_0, \dots, z_{T-1}) \coloneqq\\
     \coloneqq  \frac{1}{1+ \phi(z_{T-1}) + \dots +\phi(z_{T-1}) \dots \phi(z_{0})} \Psi(z_{T-2}) \dots \Psi(z_0) m_0, \quad \text{ for }  k = T - 1.
     \end{cases}
 \end{equation*}
\begin{exmp}
\label{ex:huang_example} A simple example is proposed in Section 3.3 of \citep{huang2013mean}. Let $T$  be equal to 2 and $(z_0,z_1)$ be given. The solution is defined by 
\[ D_0 = \frac{\phi(z_1)\phi(z_0)}{1+\phi(z_1)+\phi(z_1)\phi(z_0)}, \quad D_1 = \frac{\phi(z_1)}{1+\phi(z_1)}, \quad D_2 = 1, \]
with controls
\[ \hat{\alpha}_0(x) = \frac{(1+\phi(z_1))x}{1+\phi(z_1)+\phi(z_1)\phi(z_0)}, \quad \hat{\alpha}_1(x) = \frac{x}{1+\phi(z_1)+\phi(z_1)\phi(z_0)}, \quad \hat{\alpha}_2(x) = 0.  \]
\end{exmp}

\subsection{Solution of the MFC} \label{sec:hara_mfc_deep_solver}

We now turn our attention to the cooperative setting. For this problem, we are not aware of any explicit solution for the social optimum. Instead, we employ the numerical method proposed in~\citep{carmona2019convergence-II} and use the result as a benchmark. We recall how this method works in our context. The initial problem is to minimize over $\alpha$:
\begin{equation*}
    J(\balpha) = \mathbb{E} \sum_{t=0}^T \rho^t v(c_t) = \mathbb{E} \sum_{t=0}^T \rho^t v(X_t^{\balpha}-\alpha_t),
\end{equation*}
subject to: $X_0^{\balpha}$ has a fixed distribution and 
\begin{equation*}
    X_{t+1}^{\balpha} = G(\mathbb{E}[\alpha_t], W_t) \alpha_t, \quad t > 0.
\end{equation*}
This problem is approximated by the following one. We fix an architecture of neural network with input in $\mathbb{R}^2$ and output in $\mathbb{R}$. Such neural networks are going to play the role of the control function, in a Markovian feedback form. The inputs are the time and space variables, and the output is the value of the control. Then the goal is to minimize over parameters $\omega$ of neural networks with this architecture the following function:
\begin{equation*}
    \widetilde{J}^N(\omega) 
    = \mathbb{E}\left[\frac{1}{N}\sum_{i=1}^N \sum_{t=0}^T \rho^t v(c^i_t) \right]
    = \mathbb{E} \left[ \frac{1}{N}\sum_{i=1}^N  \sum_{t=0}^T \rho^t v(X_t^{i,\varphi_\omega}-\varphi_\omega(t,X_t^{i,\varphi_\omega})) \right],
\end{equation*}
subject to: $X_0^{i,\varphi_\omega},i=1,\dots,N$ are i.i.d. with fixed distribution and 
\begin{equation*}
    X_t^{i, \varphi_\omega} = G\left(\frac{1}{N}\sum_{j=1}^N \varphi_\omega(t,X_t^{j,\varphi_\omega}), W^i_t\right) \varphi_\omega(t,X_t^{i,\varphi_\omega}), \quad t > 0, i=1,\dots,N.
\end{equation*}
Notice that the parameters $\omega$ are used to compute the $X_t^{i, \varphi_\omega}$ for every $i$ and every $t$. The mean of the control $\mathbb{E}[\alpha_t]$ is replaced by an empirical average over $N$ samples. For this problem, an approximate minimizer is computed by running stochastic gradient descent (SGD for short) or one of its variants. At iteration $k$, we have a candidate $\omega_k$ for the parameters of the neural network. We randomly pick initial positions $\underline{X}_0 := (X_0^{i,\varphi_{\omega_k}})_{i=1,\dots,N}$ and noises $\underline{\bW} := (W^i_t)_{t=1,\dots,T, i=1,\dots,N}$. Based on this, we simulate trajectories $(X_t^{i, \varphi_{\omega_k}})_{t,i}$ and compute the associated cost, namely the term inside the expectation in the definition of $\widetilde{J}^N(\omega)$:
$$
    L(\omega_k; \underline{X}_0, \underline{\bW})
    := 
    \frac{1}{N}\sum_{i=1}^N  \sum_{t=0}^T \rho^t v(X_t^{i,\varphi_{\omega_k}}-\varphi_{\omega_k}(t,X_t^{i,\varphi_{\omega_k}})).
$$
Using backpropagation, the gradient $\nabla_\omega L(\omega_k; \underline{X}_0, \underline{\bW})$ of this cost with respect to $\omega$ is computed, and it is used to update the parameters. We thus obtain $\omega_{k+1}$ defined by:
$$
    \omega_{k+1} = \omega_{k} - \eta_k \nabla_{\omega} L(\omega_k; \underline{X}_0, \underline{\bW}),
$$
where $\eta_k>0$ is the learning rate used at iteration $k$. In our implementation for the numerical results presented below, instead of the plain SGD algorithm we used Adam optimizer~\citep{kingma2014adam}.

\subsection{Numerical results}\label{sec:results_hara}

In this section, numerical results of the application of the U2-MF-QL-FH algorithm to the mean field capital accumulation problem are presented. The interaction with the population is through the law of the controls. The algorithm \ref{algo:U2MFQL} was adapted to this case as discussed in Section \ref{sec:algorithm}.\\ The problem analyzed is a specific case of the Example \ref{ex:huang_example}. For more details we refer to \cite[Sections 6.3 and 7, Example 18]{huang2013mean}. \\
The production function is defined as follows
\begin{equation}
    G(z,W)=g(z)W, \quad \quad g(z) = \frac{1}{\rho \E[W^{\gamma}]} \frac{C}{1+(C-1)z^3},
\end{equation}
where $W$ has support $D_W =\{ 0.9, 1.3 \}$ with corresponding probabilities $[0.75, 0.25]$, $C$ is equal to 3, the discount factor $\rho$ is equal to 0.95 and the parameter $\gamma$  of the utility function defined in equation (\ref{eq:hara_utility}) is equal to 0.2. The distribution of $X_0^{\alpha,\theta}$ is uniform in $[0,1]$.\\
This problem is characterized by discrete time and continuous state and action spaces. In order to apply the U2-MF-QL-FH algorithm, these spaces are truncated and discretized as discussed in Section \ref{sec:continuous_extension_algo}. They have been chosen large enough to make sure that the state is within the boundary most of the time. In practice, this would have to be calibrated in a model-free way through experiments. In this example, for the numerical experiments, we used the knowledge of the model. \\
The action space is given by
$ \mc{A} = \{a_0=0, \dots, a_{|\mc{A}|-1}=4 \}$  and the state space by $ \mc{X} = \{x_0=0, \dots, x_{|\mc{X}|-1}=4 \}$. The step size for the discretization of the state and action spaces is
given by $0.05$. \\
The algorithm \ref{algo:U2MFQL} is adapted to this particular example. Since borrowing is not allowed, the set of admissible action at state $x$ is given by $\mc{A}(x) = \{ a\in \mc{A} \text{ if } a\leq x \} \subseteq \mc{A}$. The exploitation-exploration trade off is tackled on each episode using an $\epsilon-$greedy policy. Supposed that the agent is in state $x$, the algorithm chooses a  random action in $\mc{A}(x)$ with probability $\epsilon$  and the action in $\mc{A}(x)$ which results optimal based on the current estimation with probability $1-\epsilon$. In our example, the value of epsilon is fixed to $0.15$.  \\
\vskip 6pt
\noindent
The following numerical results show how the U2-MF-QL-FH algorithm is able to learn an approximation of the control function and the mean field term in the MFG and MFC cases depending on the choice of the parameters $(\omega^Q,\omega^{\theta})$.

 \subsubsection{Learning of the controls}
\begin{figure}[H]
\centering
\begin{minipage}{.45\textwidth} 
  \centering 
  \includegraphics[width=.9\linewidth]{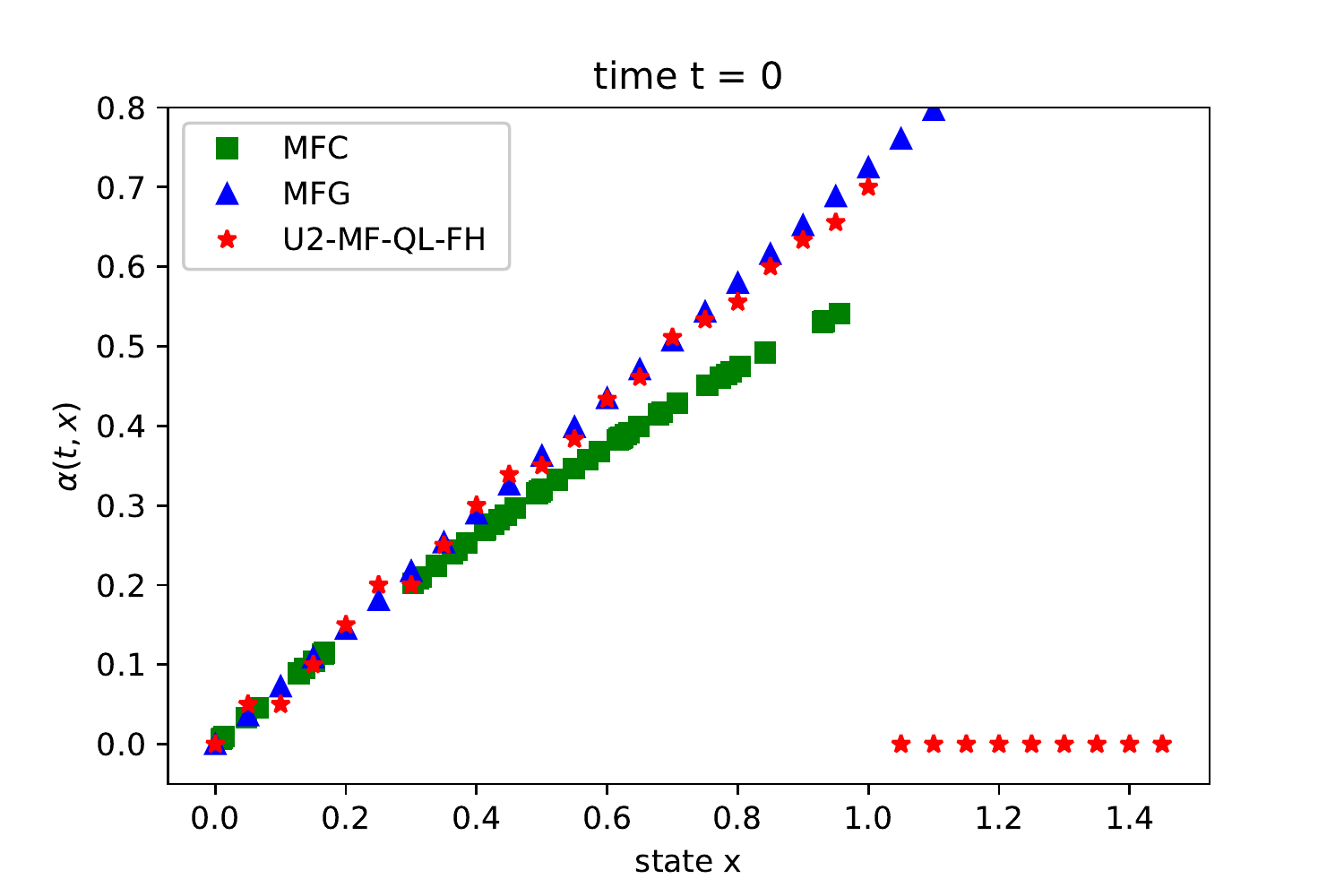}
  \caption{Learned Controls for MFG at time $0$.}
  \label{fig:control_MFG_t0}
\end{minipage}%
\hspace{0.2cm}\begin{minipage}{.45\textwidth} 
  \centering 
  \includegraphics[width=.9\linewidth]{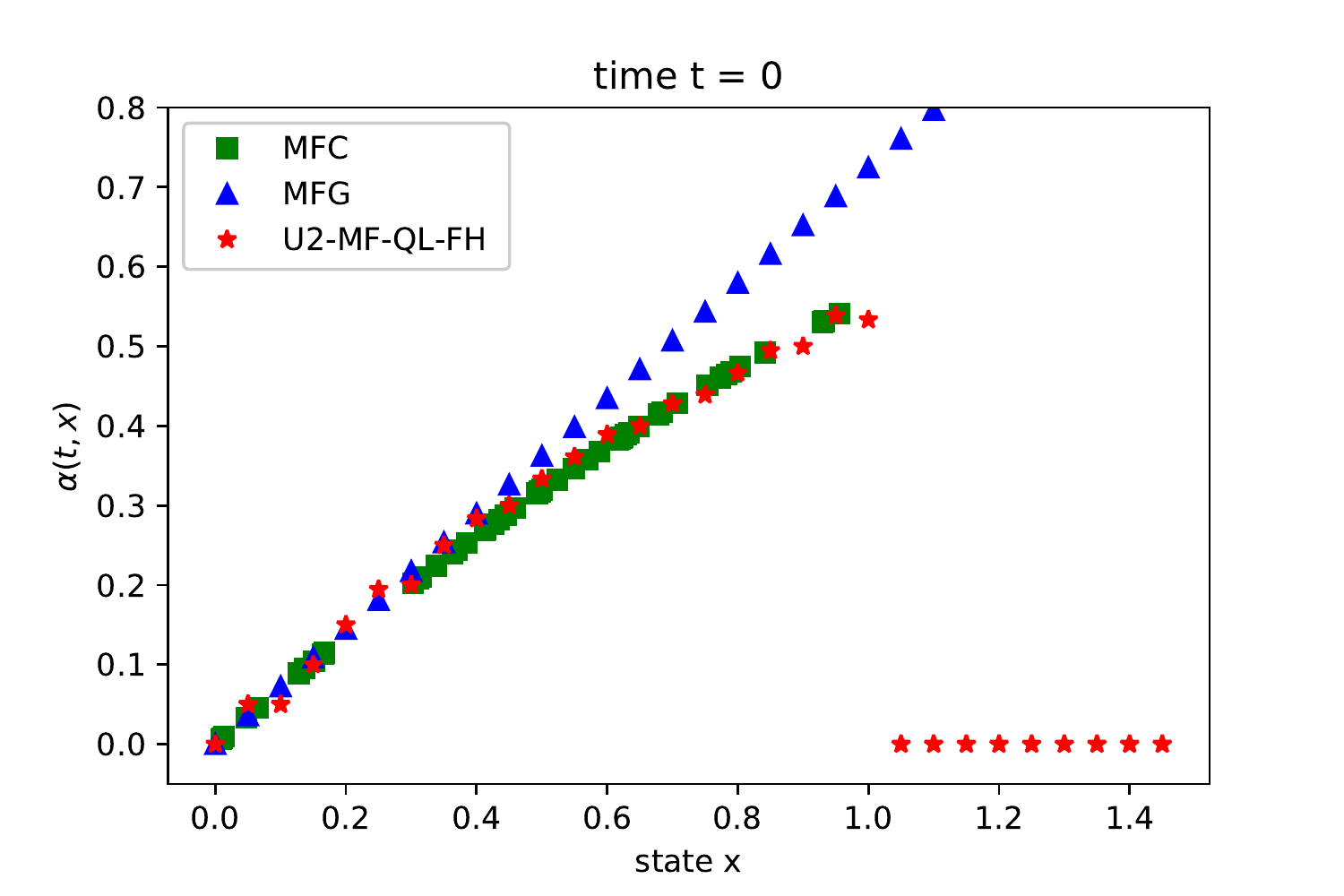}
  \caption{Learned Controls for MFC at time $0$. }
  \label{fig:control_MFC_t0}
\end{minipage}%
\end{figure}

\begin{figure}[H]
\centering
\begin{minipage}{.45\textwidth} 
  \centering 
  \includegraphics[width=.9\linewidth]{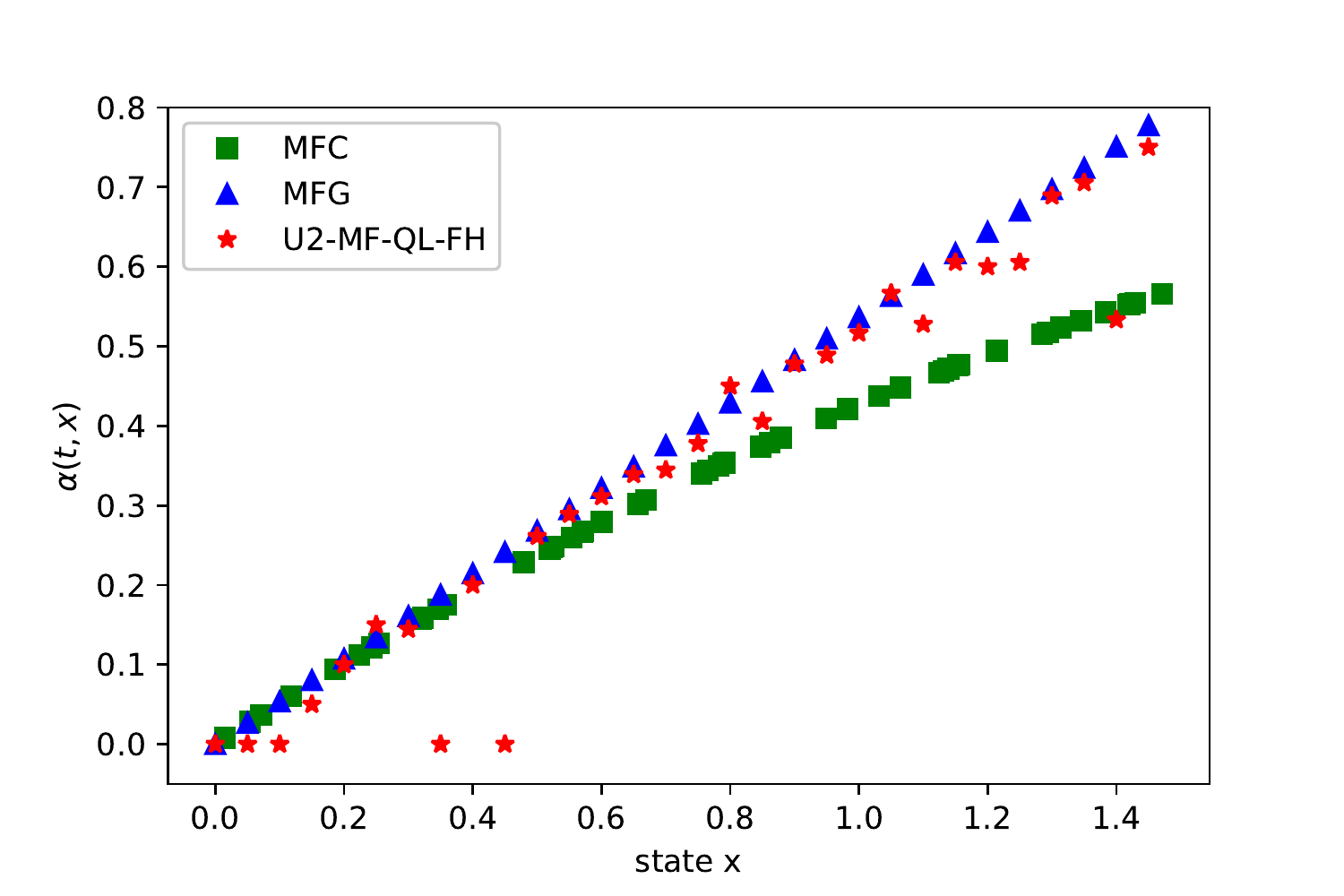}
  \caption{Learned Controls for MFG at time $1$.
  }
  \label{fig:control_MFG_t1}
\end{minipage}%
\hspace{0.2cm}\begin{minipage}{.45\textwidth} 
  \centering 
  \includegraphics[width=.9\linewidth]{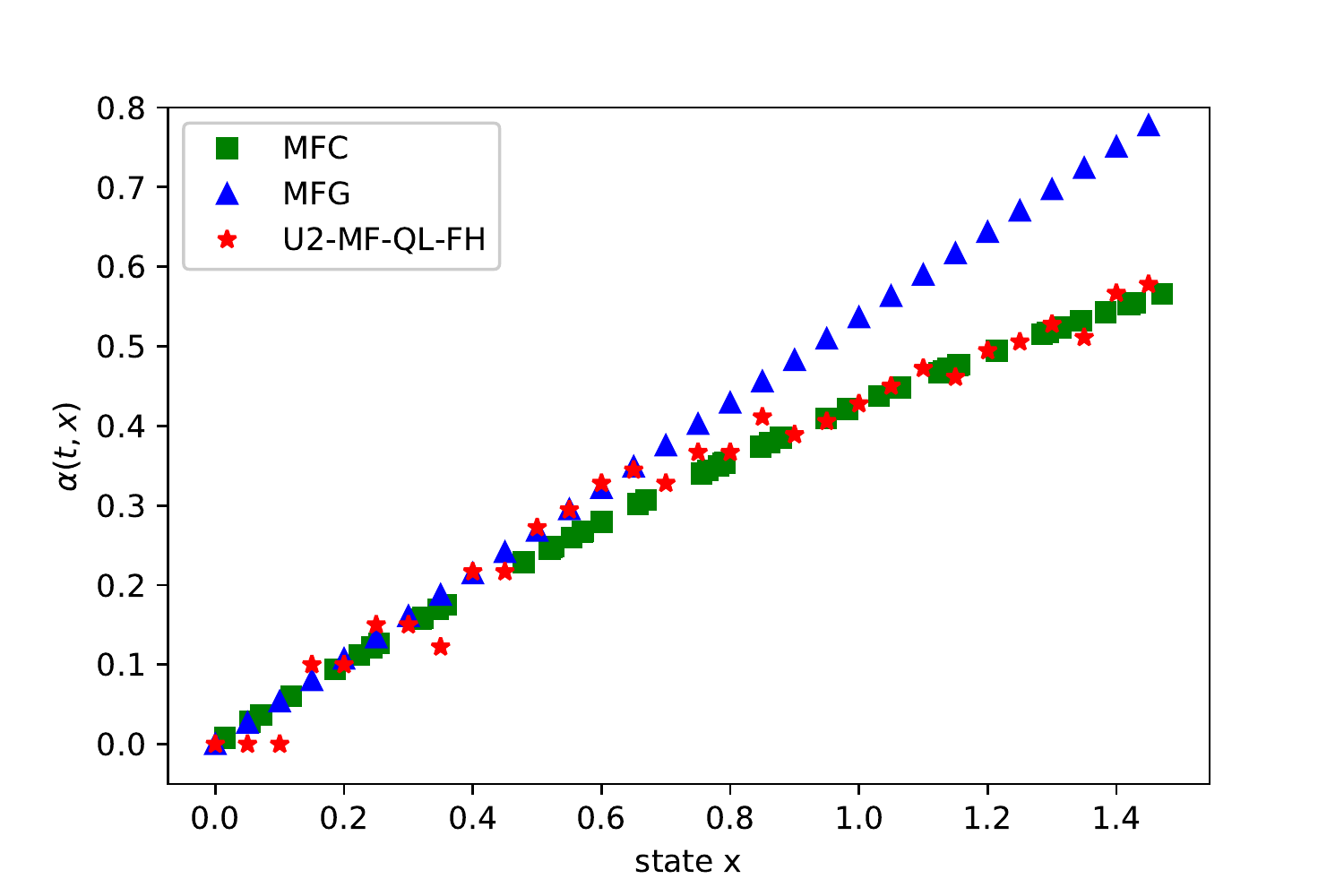}
  \caption{Learned Controls for MFC at time $1$. 
  }
  \label{fig:control_MFC_t1}
\end{minipage}%
\end{figure}

\begin{figure}[H]
\centering
\begin{minipage}{.45\textwidth} 
  \centering 
  \includegraphics[width=.9\linewidth]{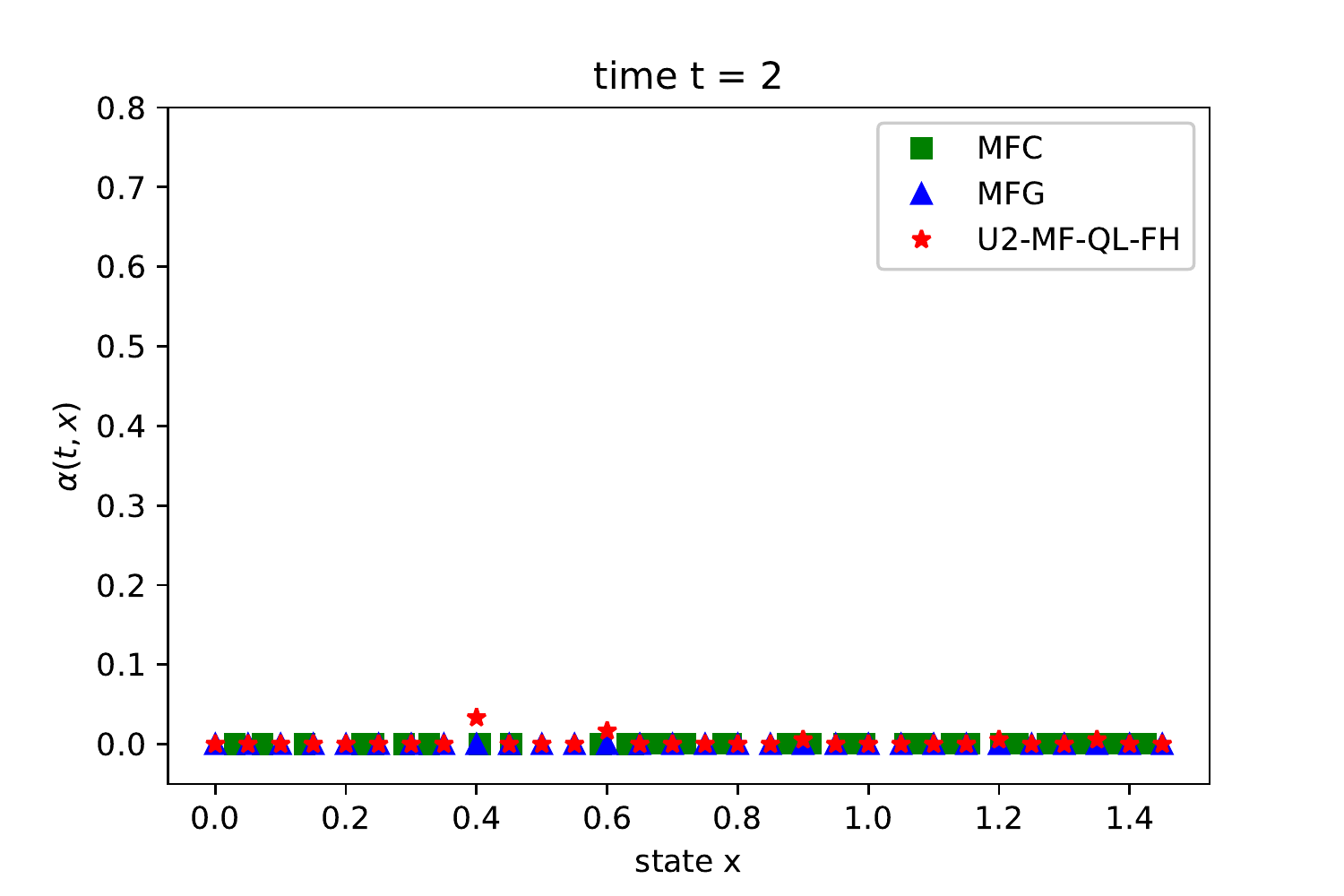}
  \caption{Learned Controls for MFG at time $2$.
  }
  \label{fig:control_MFG_t2}
\end{minipage}%
\hspace{0.2cm}\begin{minipage}{.45\textwidth} 
  \centering 
  \includegraphics[width=.9\linewidth]{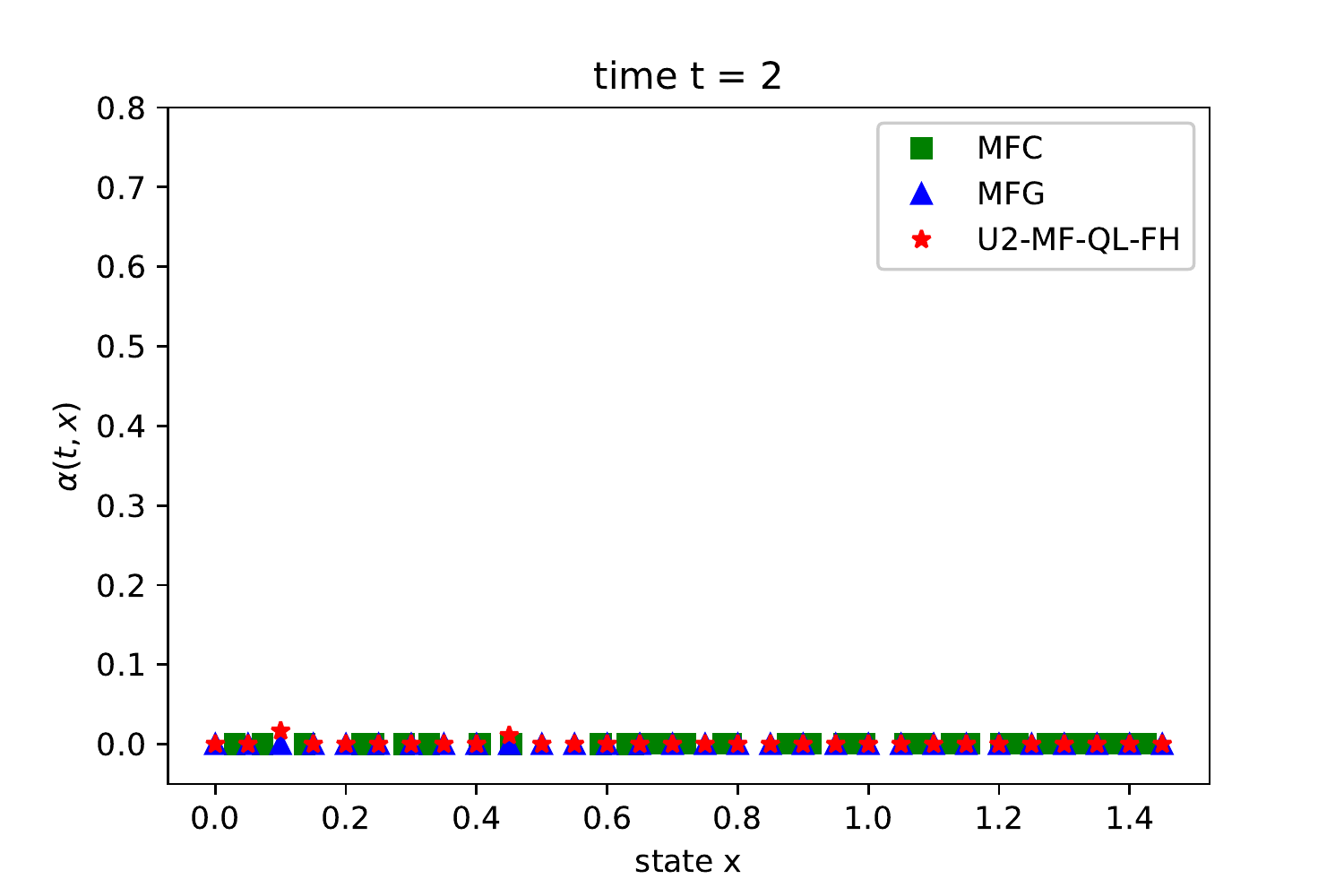}
  \caption{Learned Controls for MFC at time $2$.
  }
  \label{fig:control_MFC_t2}
\end{minipage}%
\end{figure}

\textbf{Figures \ref{fig:control_MFG_t0}, \ref{fig:control_MFC_t0}, \ref{fig:control_MFG_t1}, \ref{fig:control_MFC_t1},\ref{fig:control_MFG_t2}, \ref{fig:control_MFC_t2}:  controls learned by the algorithm.} The controls learned by the U2-MF-QL-FH algorithm are compared with the benchmark solutions. Each plot corresponds to a different time point $t \in \{0,1,2\}$. The $x-$axis represents the state variable $x$. The $y-$axis relates to the action $\alpha_t(x)$. The blue (resp. green) markers show the benchmark control function for the MFG (resp. MFC) problem. The red markers are the controls learned by the algorithm. The plots show how the algorithm converges to different solutions based on the choice of the pair $(\omega^Q,\omega^{\theta})$. On the left, the choice $(\omega^Q,\omega^{\theta})=(0.55,0.85)$ produces the approximation of the solution of the MFG. On the right, the set of parameters  $(\omega^Q,\omega^{\theta})=(0.7,0.05)$ lets the algorithm learn the solution of the MFC problem. The results presented in the Figures are averaged over 10 runs.

\subsubsection{Learning of the mean field}
\begin{figure}[H]
\centering
\begin{minipage}{.45\textwidth} 
  \centering 
  \includegraphics[width=.9\linewidth]{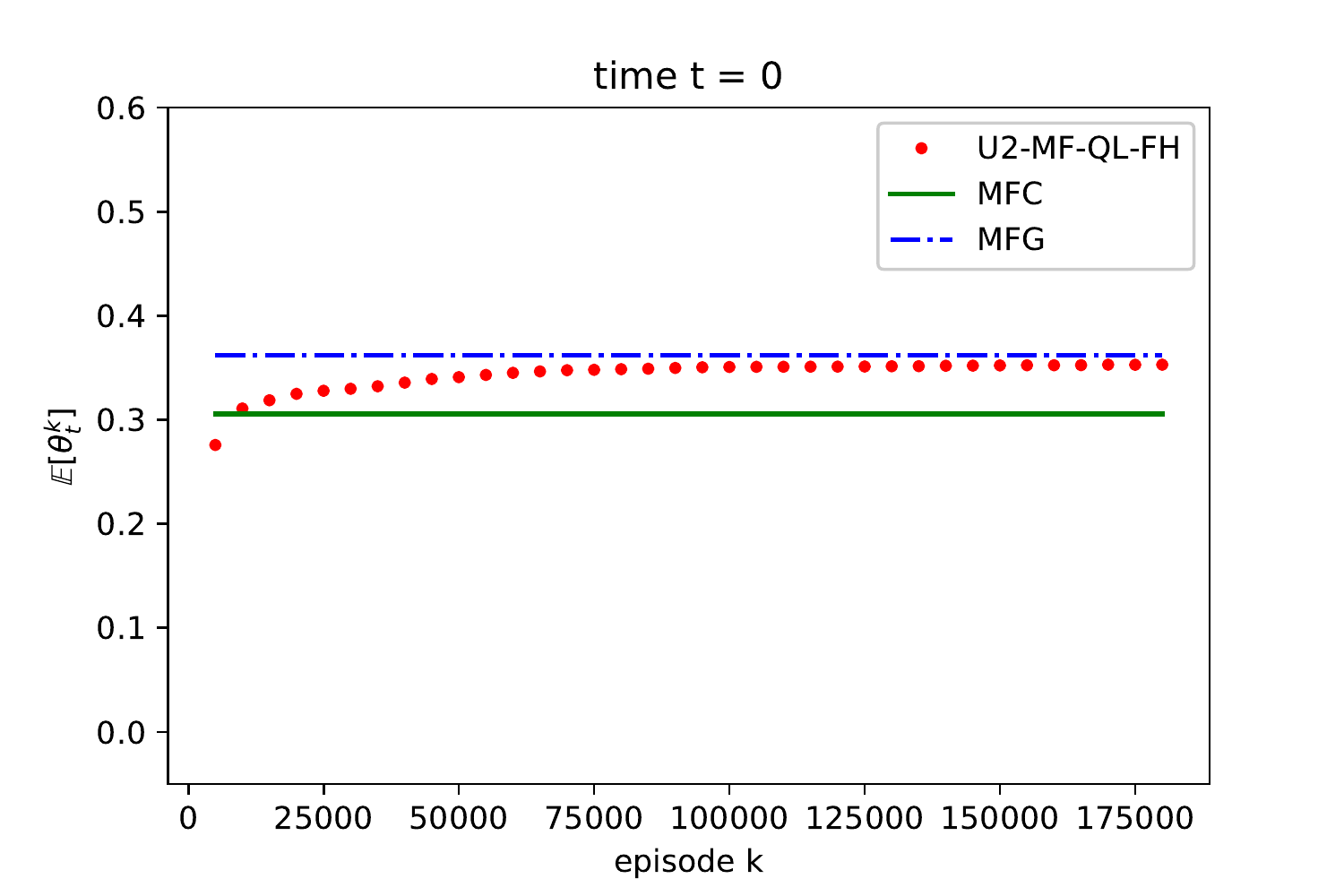}
  \caption{Learned control's mean for MFG at time $0$.
  }
  \label{fig:theta_MFG_t0}
\end{minipage}%
\hspace{0.2cm}\begin{minipage}{.45\textwidth} 
  \centering 
  \includegraphics[width=.9\linewidth]{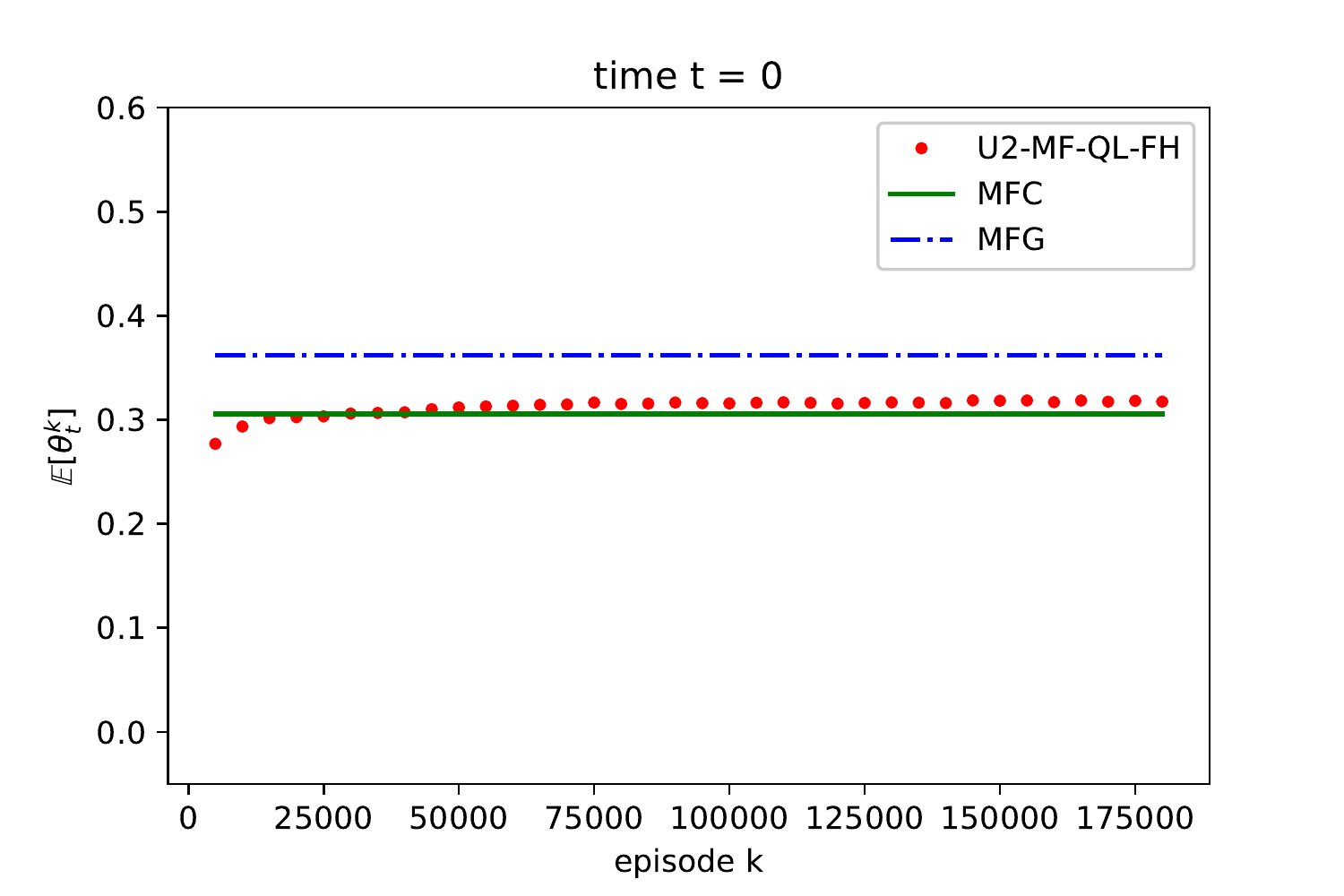}
  \caption{Learned control's mean for MFC at time $0$.
  }
  \label{fig:theta_MFC_t0}
\end{minipage}%
\end{figure}

\begin{figure}[H]
\centering
\begin{minipage}{.45\textwidth} 
  \centering 
  \includegraphics[width=.9\linewidth]{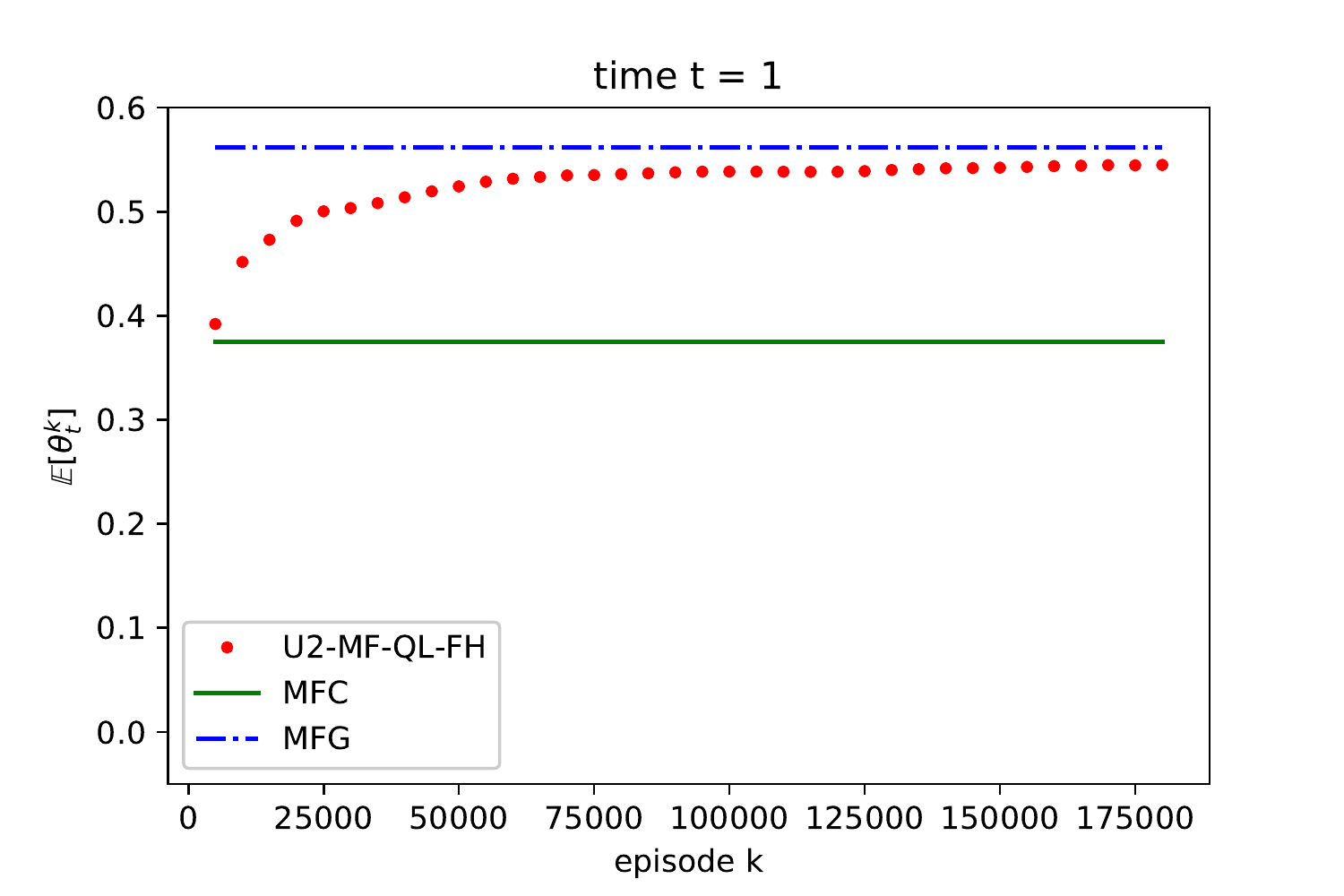}
  \caption{Learned control's mean for MFG at time $1$.
  }
  \label{fig:theta_MFG_t1}
\end{minipage}%
\hspace{0.2cm}\begin{minipage}{.45\textwidth} 
  \centering 
  \includegraphics[width=.9\linewidth]{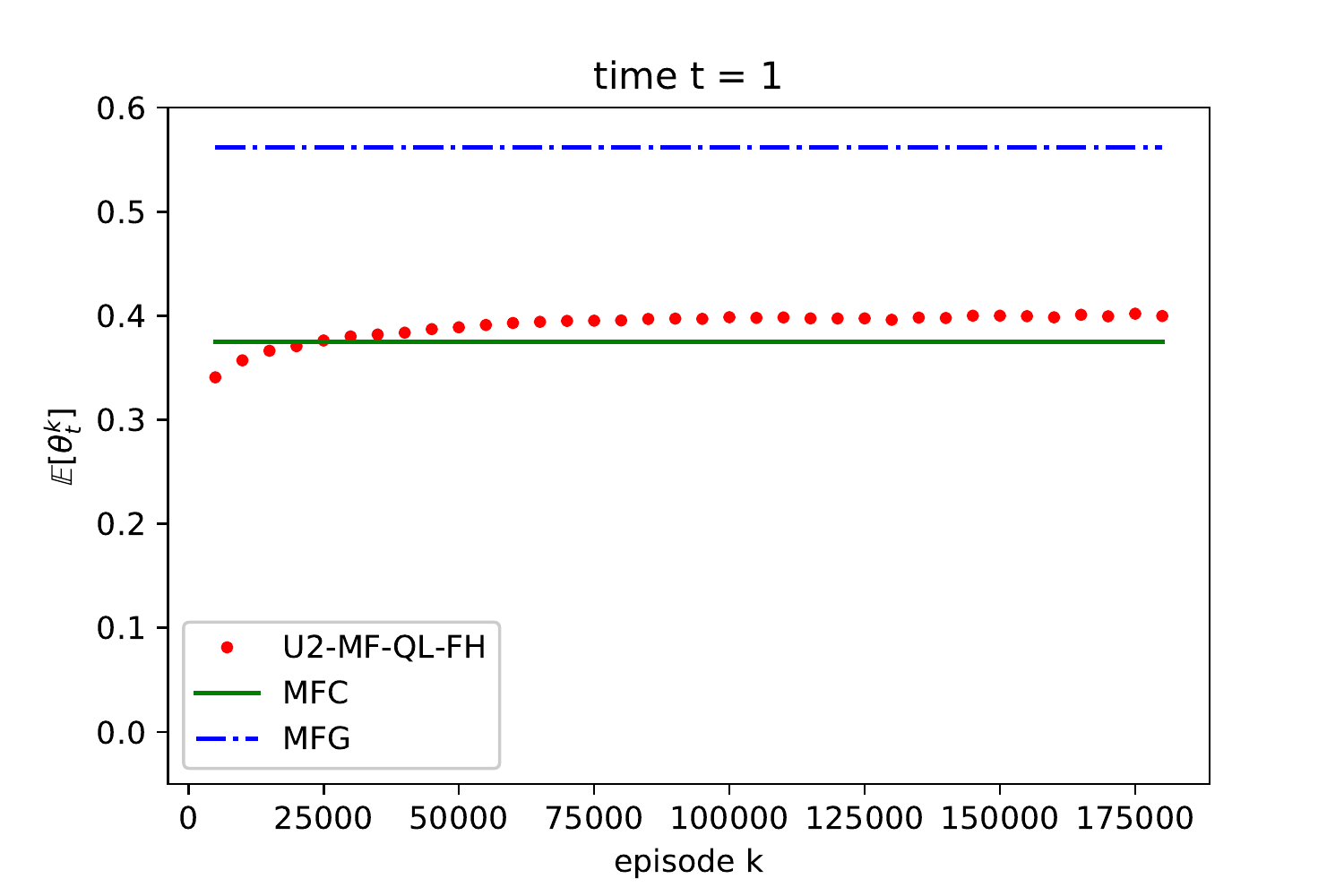}
  \caption{Learned control's mean for MFC at time $1$. 
  }
  \label{fig:theta_MFC_t1}
\end{minipage}%
\end{figure}

\begin{figure}[H]
\centering
\begin{minipage}{.45\textwidth} 
  \centering 
  \includegraphics[width=.9\linewidth]{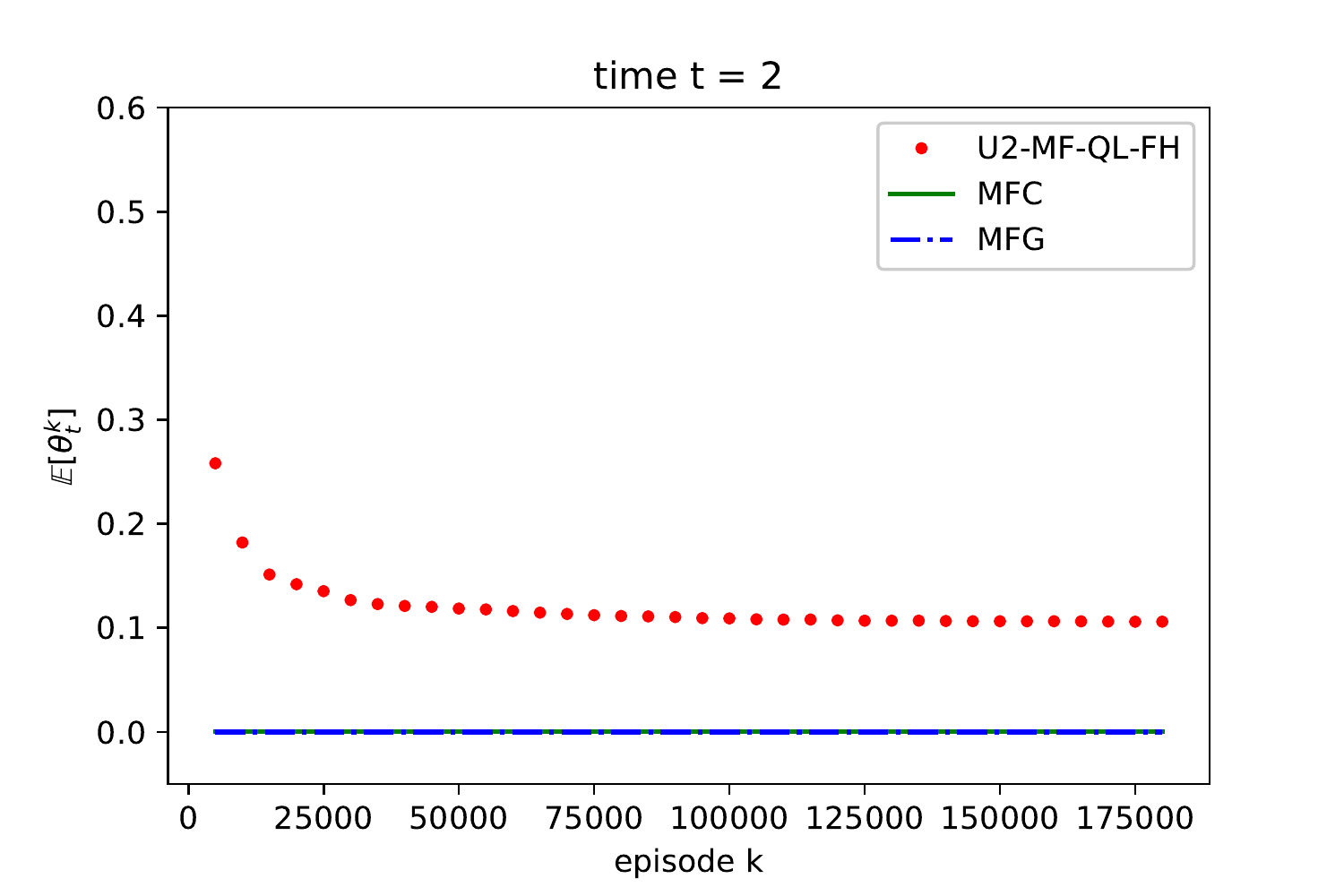}
  \caption{Learned control's mean for MFG at time $2$.
  }
  \label{fig:theta_MFG_t2}
\end{minipage}%
\hspace{0.2cm}\begin{minipage}{.45\textwidth} 
  \centering 
  \includegraphics[width=.9\linewidth]{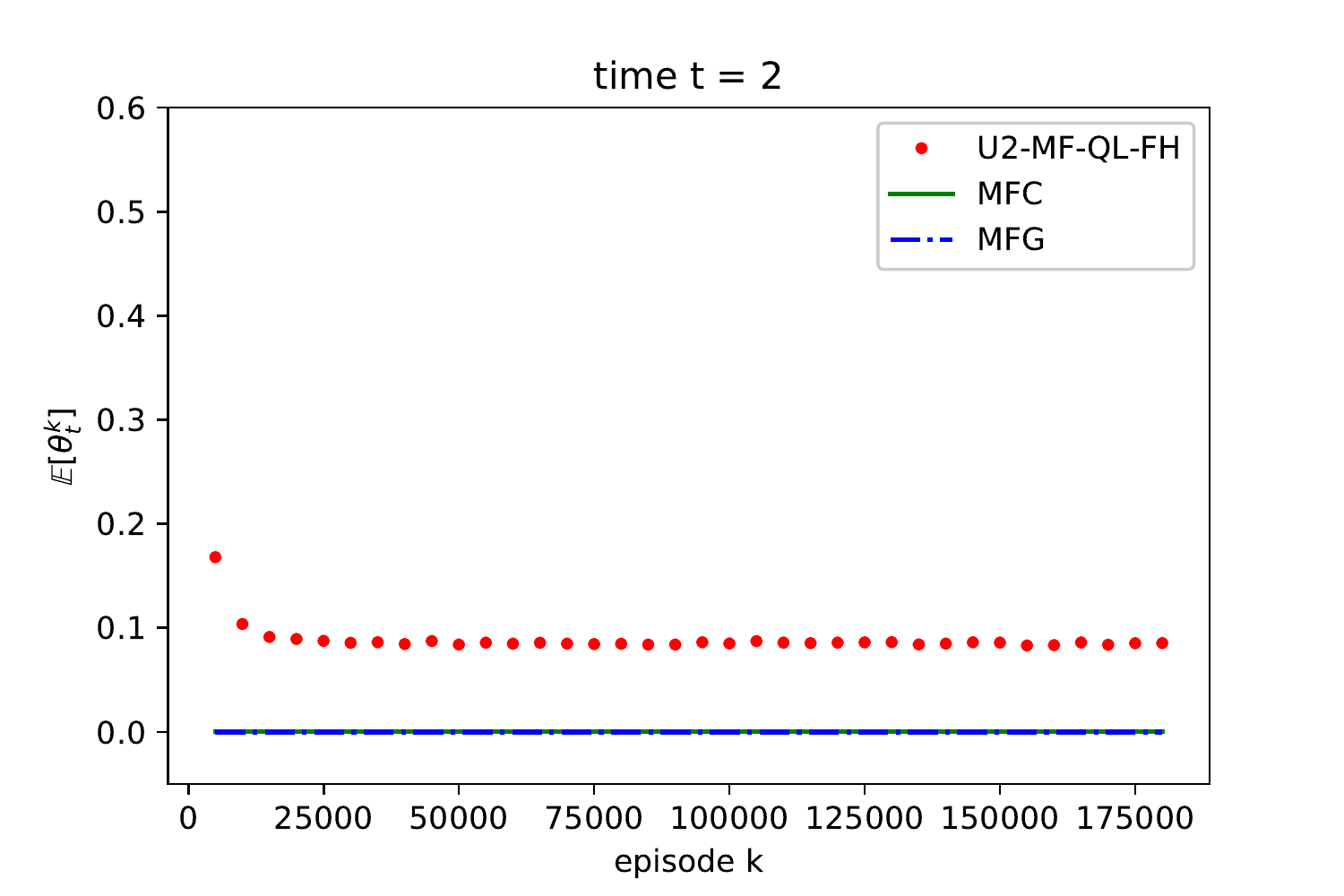}
  \caption{Learned control's mean for MFC at time $2$. 
  }
  \label{fig:theta_MFC_t2}
\end{minipage}%
\end{figure}

\textbf{Figures \ref{fig:theta_MFG_t0}, \ref{fig:theta_MFC_t0}, \ref{fig:theta_MFG_t1}, \ref{fig:theta_MFC_t1}, \ref{fig:theta_MFG_t2}, \ref{fig:theta_MFC_t2}:  $\E[\alpha_t]$ learned by the algorithm.} 
The estimation of the first moment of the distribution of the controls evolves with respect the number of learning episodes. The estimated quantity is compared with the benchmarks presented in Sections \ref{sec:hara_mfg_sol} and \ref{sec:hara_mfc_deep_solver}. Each plot corresponds to a different time point $t \in \{0,1,2\}$. The $x-$axis represents the learning episode $k$. The $y-$axis relates to the estimate of the first moment of the mean field $\E[\alpha^k_t]$ obtained by episode k.  The blue (resp. green) line shows the benchmark solution for the MFG (resp. MFC) problem. The red dots are the estimates learned by the algorithm. On the left, the algorithm reaches the solution of the MFG based on the parameters $(\omega^Q,\omega^{\theta})=(0.55,0.85)$. On the right, the values  $(\omega^Q,\omega^{\theta})=(0.7,0.05)$ allows the algorithm to converge to the solution of the MFC problem. The results presented in the Figures are averaged over 10 runs.

\section{A mean field execution problem}\label{sec:trader}

We now consider the \textit{Price Impact Model} as an example of application to finance originally studied by Carmona and Lacker in \citep{MR3325272}, and presented in the book of Carmona and Delarue  \cite[Sections 1.3.2 and 4.7.1]{MR3752669}. This model addresses the question of optimal execution in the context of high frequency trading when a large group of traders want to buy or sell shares before a given time horizon $T$ (\textit{e.g.}, one day). The price of the stock is influenced by the actions of the traders: if they buy, the price goes up, whereas if they sell, the price goes down. This effect is stronger if a significant proportion of traders buy or sell at the same time. Incorporating such a price impact naturally leads to a problem with mean field interactions through the traders' actions.

 Approaching this problem as a mean field game, the inventory of the representative trader is modeled by a stochastic process $(X_t)_{0 \leq t \leq T} $ such that
\begin{equation*}
dX_t = \alpha_t dt +\sigma dW_t, \quad t\in [0,T],
\end{equation*}
where $\alpha_t$ corresponds to the trading rate and $W$ is a standard Brownian motion. The price of the asset $(S_t)_{0 \leq t \leq T}$ is influenced by the trading strategies of all the traders through the mean of the law of the controls $(\theta_t=\mathcal{L}(\alpha_t))_{0 \le t \le T}$ as follows:
\begin{equation*}
dS_t =\gamma \biggl( \int_{\mathbb{R}}a d\theta_t(a) \biggr) dt + \sigma_0 dW_t^0, \quad t\in [0,T],
\end{equation*} 
where $\gamma$ and $\sigma_0$ are constants and the Brownian motion $W^0$ is independent from $W$.
The amount of cash held by the trader at time $t$ is denoted by the process $(K_t)_{0 \le t \le T}$. The dynamic of $K$ is modeled by
\begin{equation*}
dK_t=-[\alpha_t S_t +c_{\alpha}(\alpha_t)]dt,
\end{equation*}
where the function $\alpha \mapsto c_{\alpha}(\alpha)$ is a non-negative convex function satisfying $c_{\alpha}(0)=0$, representing the cost for trading at rate $\alpha$. The wealth $V_t$ of the trader at time $t$ is defined as the sum of the cash held by the trader and the value of the inventory with respect to the price $S_t$:
\begin{equation*}
V_t=K_t+X_t S_t.
\end{equation*}
Applying the self-financing condition of Black-Scholes' theory, the changes over time of the wealth $V$ are given by the equation:
\begin{equation}
\begin{split}
    dV_t&=dK_t+X_tdS_t+S_t dX_t
    \\
   & =\Big[ -c_{\alpha}(\alpha_t)+\gamma X_t \int_{\mathbb{R}} a d\theta_t(a) \Big]dt + \sigma S_t dW_t + \sigma_0 X_t dW_t^0.
\end{split}
\label{wealth}
\end{equation}
We assume that the trader is subject to a running liquidation constraint modeled by a function $c_X$ of the shares they hold, and to a terminal liquidation constraint at maturity $T$ represented by a scalar function $g$. Thus, the cost function is defined by:
\begin{equation*}
J(\alpha)=\mathbb{E}\Big[\int_0^T c_X(X_t) dt +g(X_T) - V_T\Big],
\end{equation*}
where the terminal wealth $V_T$ is taken into account with a negative sign as the cost function is to be minimized.
From equation (\ref{wealth}), it follows that
\begin{equation*}
J(\alpha)=\mathbb{E}\Big[ \int_0^T f(t,X_t,\theta_t,\alpha_t)dt +g(X_T)\Big],    
\end{equation*}
where the running cost is defined by
\begin{equation*}
f(t,x,\theta,\alpha)=c_{\alpha}(\alpha)+c_X(x)-\gamma x \int_{\mathbb{R}} a  d\theta(a),
\end{equation*}
for $0\leq t \leq T$, $x \in \mathbb{R}^d$, $\theta \in \mathcal{P}(\mathbb{A})$ and $\alpha \in \mathbb{A} = \mathbb{R}$. We assume that the functions $c_X$ and $g$ are quadratic and that the function $c_{\alpha}$ is strongly convex in the sense that its second derivative is bounded away from $0$. Such a particular case is known as the Almgren-Chriss linear price impact model. Thus, the control is chosen to minimize:
\begin{equation*}
J(\alpha)=\mathbb{E}\left[ \int_0^T \left( \frac{c_{\alpha}}{2}{\alpha_t}^2+\frac{c_X}{2}X_t^2-\gamma X_t\int_{\mathbb{R}} a d\theta_t(a) \right)dt + \frac{c_g}{2}X_T^2\right],
\end{equation*}
over $\alpha \in \mathbb{A}$.
To summarize, the running cost consists of three components. The first term represents the cost for trading at rate $\alpha$. The second term takes into consideration the running liquidation constraint in order to penalize unwanted inventories. The third term defines the actual price impact. Finally, the terminal cost represents the terminal liquidation constraint.

\subsection{The MFG trader problem}
Referring to Section \ref{sec:mfg}, 
the MFG problem is solved by first solving a standard stochastic control problem where the flow of distribution of control is given and then, solving a fixed point problem ensuring that this flow of distribution is identical to the flow of distributions of the optimal control. We adopt here the FBSDE approach where the backward variable represents the derivative of the value function. In other words, the optimal control is obtained by minimizing the Hamiltonian
\begin{equation}\label{Hamiltonian}
    H(x,\alpha,\theta,y) 
    = \left(\frac{c_{\alpha}}{2}{\alpha}^2+\frac{c_X}{2}x^2-\gamma x\int_{\mathbb{R}} a d\theta(a) \right) + \alpha y,
\end{equation}
to obtain
\begin{equation}
\label{eq:trader-MFG-ctrl}
\hat\alpha_t = - \frac{1}{c_{\alpha}}Y_t,
\end{equation}
where $(X,Y)$ solves the FBSDE system obtained via the Pontryagin approach:
\begin{equation}
\label{eq:trader-MFG-FBSDE-pontryagin}
\left\{
\begin{aligned}
    dX_t &= -\frac{1}{c_{\alpha}}Y_t dt + \sigma dW_t, \quad \quad X_0 \sim \mu_0 \\
    dY_t &= -\left(c_X X_t + \frac{\gamma }{c_{\alpha}} \mathbb{E} [Y_t]\right) dt + Z_t dW_t,\quad Y_T = c_g X_T. \\
    \end{aligned}
\right.
\end{equation}

\subsection*{Solution of the MFG problem} \label{sec:mfg-sol}

The solution of the mean field game case is discussed in details in \cite[Sections 1.3.2 and 4.7.1]{MR3752669}. 
In a nutshell, one takes expectation in \eqref{eq:trader-MFG-FBSDE-pontryagin} to obtain a system of forward-backward ODEs for the mean of $X_t$ denoted by $\bar{x}_t$ and the mean of $Y_t$ denoted by $\bar{y}_t$.
This system is solved using the ansatz $\bar{y}_t = \bar{\eta}_t \bar{x}_t + \bar{\chi}_t$. The coefficient function $\bar{\eta}_t$ satisfies a Riccati equation which admits the solution:
\begin{equation*}
\bar{\eta}_t= \frac{-C (e^{(\delta^+-\delta^-)(T-t)}-1)-c_g(\delta^+e^{(\delta^+-\delta^-)(T-t)}-\delta^-)}{(\delta^-e^{(\delta^+-\delta^-)(T-t)}-\delta^+)-c_gB(e^{(\delta^+-\delta^-)(T-t)}-1)},
\end{equation*}
for $t\in[0,T]$, where $B=1/c_{\alpha}$, $C=c_X, \delta^\pm=-D \pm \sqrt{R}$, 
with
$D = -\gamma /(2c_{\alpha})$,
$R=D^2+BC$ and $\bar{x}_0 =\E[X_0]$. Additionally, we found $\bar{\chi}_t = 0$, and 
$$\bar{x}_t = \bar{x}_0 e^{-\int_0^t \frac{\bar{\eta}_s}{c_{\alpha}}ds}.$$

 The FBSDE system (\ref{eq:trader-MFG-FBSDE-pontryagin}) is solved by replacing $\E[Y_t]$  with the explicit expression for $\bar{y}_t = \bar{\eta}_t \bar{x}_t + \bar{\chi}_t$, and  using the ansatz $Y_t=\eta_t X_t+\chi_t$. One finds the following explicit formulas for the coefficient functions $\eta_t$ and $\chi_t$:
  \begin{equation*}
\begin{split}
\eta_t&=-c_{\alpha}\sqrt{c_X/c_\alpha}\frac{c_{\alpha}\sqrt{c_X/c_\alpha}-c_g-(c_{\alpha}\sqrt{c_X/c_\alpha}+c_g)e^{2\sqrt{c_X/c_\alpha}(T-t)}}{c_{\alpha}\sqrt{c_X/c_\alpha}-c_g+(c_{\alpha}\sqrt{c_X/c_\alpha}+c_g)e^{2\sqrt{c_X/c_\alpha}(T-t)}}, \\
\chi_t &= (\bar{\eta}_t - \eta_t ) \bar{x}_t. \\
\end{split}
\end{equation*}
Finally, the optimal control (\ref{eq:trader-MFG-ctrl}) is given by
$\hat\alpha_t =\hat\alpha (t, X_t)$ where
\begin{equation}
    \hat{\alpha}(t,x) = -\frac{1}{c_{\alpha}} \left ( \eta_t x + (\bar{\eta}_t - \eta_t ) \bar{x}_t \right).
\end{equation}

\subsection{The MFC trader problem}

In the case of mean field control (\textit{i.e.}, control of McKean-Vlasov dynamics), following~\cite[Theorem 3.2]{acciaio2018extended} and~\cite[Section 5.3.2]{lauriere2020convergence}, 
we find that the optimal control is given by
\begin{equation}
\label{eq:trader-MFC-ctrl}
    \alpha^*_t = -\frac{1}{c_{\alpha}}\left(Y_t  - \gamma \mathbb{E}[X_t]\right),
\end{equation}
which differs from the equilibrium control~\eqref{eq:trader-MFG-ctrl} from the MFG solution because the optimality condition in the MFC case involves the derivative of the Hamiltonian (\ref{Hamiltonian}) with respect to the distribution of controls. More precisely, we have
$$
    0 
    = \partial_\alpha  H(X_t,\alpha_t,\theta_t,Y_t) + \tilde{\mathbb{E}}\left[ \partial_\theta  H(\tilde X_t, \tilde \alpha_t, \tilde \theta_t, \tilde Y_t)(\alpha_t)\right]
    = c_\alpha \alpha_t + Y_t - \gamma \mathbb{E}[X_t].
$$
Then, the corresponding FBSDE system becomes
\begin{equation}
\label{eq:trader-MFC-FBSDE-pontryagin}
\left\{
\begin{aligned}
    dX_t &= - \frac{1}{c_{\alpha}} \left(Y_t  - \gamma \mathbb{E}[X_t]\right)dt + \sigma dW_t, \quad X_0 \sim \mu_0 \\
    dY_t &= -\left(c_X X_t + \frac{\gamma }{c_{\alpha}} \mathbb{E} [Y_t] - \frac{\gamma^2}{c_\alpha} \mathbb{E}[X_t]\right) dt + Z_t dW_t,\quad Y_T = c_g X_T. \\
    \end{aligned}
\right.
\end{equation}
As a consequence, the two FBSDE systems~\eqref{eq:trader-MFG-FBSDE-pontryagin} and~\eqref{eq:trader-MFC-FBSDE-pontryagin} respectively for MFG and MFC differ.

\subsection*{Solution of the MFC problem}
The approach to obtain the solution of the MFC problem is similar to what was presented in Section \ref{sec:mfg-sol} for the MFG, but taking into consideration the extra terms due to the derivative of the Hamiltonian with respect to the distribution of controls.

First, taking expectation in \eqref{eq:trader-MFC-FBSDE-pontryagin}, one obtains the following system of forward-backward ODEs:
\begin{equation}
\label{eq:mean-trader-MFC-FBSDE-pontryagin}
\left\{
\begin{aligned}
    \dot{\bar{x}}_t &= - \frac{1}{c_{\alpha}} \left(\bar{y}_t  - \gamma \bar{x}_t\right) , \quad \bar{x}_0=x_0, \\
    \dot{\bar{y}}_t &= -\left(c_X \bar{x}_t + \frac{\gamma }{c_{\alpha}} \bar{y}_t - \frac{\gamma^2}{c_\alpha} \bar{x}_t\right)  ,\quad \bar{y}_T = c_g \bar{x}_T. \\
    \end{aligned}
\right.
\end{equation}
Using the ansatz  $\bar{y}_t = \bar{\phi}_t \bar{x}_t + \bar{\psi}_t$, we deduce that the coefficient functions $\bar\phi_t$ and $\bar\psi_t$ must satisfy
\begin{equation}
\label{eq:mean2-trader-MFC-FBSDE-pontryagin}
\left\{
\begin{aligned}
    &\dot{\bar{\phi}}_t  + 2 \frac{\gamma}{c_{\alpha}} \bar{\phi}_t - \frac{1}{c_{\alpha}} \bar{\phi}_t^2 + c_X -\frac{\gamma^2}{c_{\alpha}} , \quad \bar{\phi}_T=c_g, \\
    &\dot{\bar{\psi}}_t + \frac{1}{c_{\alpha}}(\gamma - \bar{\phi}_t) \bar{\psi}_t = 0 ,\quad \bar{\psi}_T = 0. \\
    \end{aligned}
\right.
\end{equation}
From the second equation we get $\bar\psi_t=0$ for all $t\in[0,T]$, and solving the Riccati equation for $\bar\phi_t$, we obtain:
\begin{equation} \label{eq:mfc_eta_bar}
\bar{\phi}_t
=-\frac{1}{R}\frac{(c_2 + R c_g)c_1 e^{(T-t)(c_2-c_1)}- c_2(c_1 + R c_g)}{ (c_2 + R c_g) e^{(T-t)(c_2-c_1)}- (c_1 + R c_g)},
\end{equation}
where $c_{1/2}=\frac{-a \pm \sqrt{a^2 - 4 b}}{2}$ are the roots of
$c^2 + a c + b = 0$, with $a=2\gamma R,\,b=R(\gamma^2 R-c_X)$, and $R=1/c_\alpha$.

Using $\bar{y}_t = \bar{\phi}_t \bar{x}_t$ in the first equation of (\ref{eq:mean-trader-MFC-FBSDE-pontryagin}), we obtain a first-order linear equation for $\bar{x}_t$ which admits the solution
\begin{equation} \label{eq:mfc_soln_x_bar}
    \bar x_t = \bar{x}_0 e^{-\frac{1}{c_{\alpha}} \left( \int_0^t \bar \phi_s ds -\gamma t \right)}.
\end{equation}

The solution of the McKean-Vlasov FBSDE system (\ref{eq:trader-MFC-FBSDE-pontryagin}) is obtained using the ansatz $Y_t = \phi_t X_t + \psi_t$. Observe that the drift terms in the equations for $Y_t$ in the systems (\ref{eq:trader-MFG-FBSDE-pontryagin}) and (\ref{eq:trader-MFC-FBSDE-pontryagin}) have the same linear component $-c_X X_t$. Due to this similarity, the slope coefficient functions $\eta_t$ and $\phi_t$ are identical;
$$\eta_t=\phi_t,\quad \mbox{for all}\quad t\in[0,T].$$
However, the function $\psi_t = (\bar{\phi}_t - \phi_t ) \bar{x}_t $ differs  from $\chi_t$ in the MFG case due to the new formulations of $\bar{\phi}_t$ and $\bar{x}_t$ given in (\ref{eq:mfc_eta_bar}) and (\ref{eq:mfc_soln_x_bar}). Finally, the optimal control (\ref{eq:trader-MFC-ctrl}) is given by $\alpha^*_t=\alpha^*(t,X_t)$ where
\begin{equation}
   {\alpha}^*(t,x) = -\frac{1}{c_{\alpha}} \left ( \phi_t x + (\bar{\phi}_t - \phi_t - \gamma) \bar{x}_t \right).
\end{equation}

\subsection{Numerical results}

In this section, numerical results of the application of the U2-MF-QL-FH algorithm to the trader problem are discussed. As in the case of the mean field capital accumulation problem, the interaction with the population is through the law of the controls. The algorithm \ref{algo:U2MFQL} is adapted to this case as discussed in Section \ref{sec:algorithm}. \\
We consider the problem defined by the choice of parameters: $c_{\alpha}=1$, $c_x=2$, $\gamma=1.75$, and $c_g=0.3$. The time horizon is equal to $T=1$. The distribution of the inventory process at initial time $X_{t_0}$ is Gaussian with mean $0.5$ and standard deviation $0.3$. The volatility of the process $X_t$ is given by $\sigma=0.5$.  \\
This problem is characterized by continuous time and continuous state and action spaces. In order to solve this problem using the U2-MF-QL-FH algorithm, truncation and discretization techniques together with a projector operator are applied. The time interval $[0,T]$ is uniformly discretized as $\tau = \{ t_0, \dots, t_{N_T} = T\}$ with  $\Delta t = 1/16$. The state and action spaces are truncated and discretized as discussed in Section \ref{sec:continuous_extension_algo}. The truncation parameters are chosen large enough to make sure that the state is within the boundary most of the time.  \\
In the MFG (resp. MFC), the action space is given by $\mc{A} = \{ a_0=-2.5, \dots, a_{|\mc{A}|-1}=1 \}$ (resp. $\mc{A} = \{ a_0=-0.25, \dots, a_{|\mc{A}|-1}=5 \}$) and the state space by $\mc{X} = \{ x_0=-1.5, \dots, x_{|\mc{X}|-1}=1.75 \}$ (resp. $\mc{X} = \{ x_0=-0.75, \dots, x_{|\mc{X}|-1}=4 \}$). The step size for the discretization of the spaces $\mc{A}$, and $\mc{X}$ is given by $\Delta_{a}=\Delta_{x} = \sqrt{\Delta t} = 1/4 $. The exploitation-exploration trade off is tackled on each episode  using an $\epsilon-$greedy policy. Supposed the agent is in state $x$, the algorithm picks the action that is optimal based on the current estimates with probability $1-\epsilon$ and a random action in $\mc{A}$ with probability $\epsilon$. In particular, the value of epsilon is fixed to $0.1$.

\vskip 6pt
\noindent
The following numerical results show how the U2-MF-QL-FH algorithm is able to learn an approximation of the control function and the mean field term in the MFG and MFC cases depending on the choice of the parameters $(\omega^Q,\omega^{\theta})$.

 \subsubsection{Learning of the controls}
\begin{figure}[H]
\centering
\begin{minipage}{.45\textwidth} 
  \centering 
  \includegraphics[width=.9\linewidth]{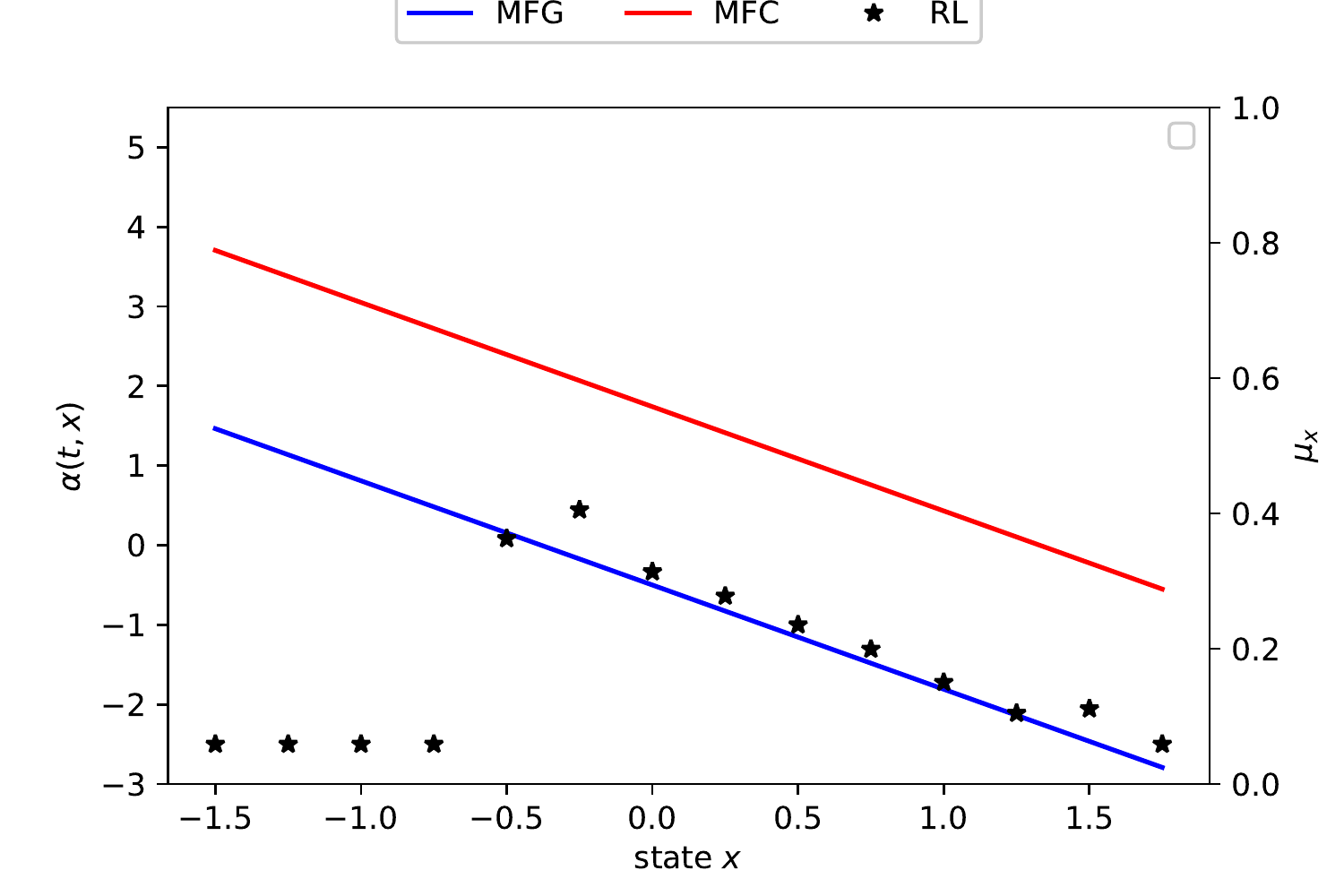}
  \caption{Learned Controls for MFG at time $0$.
  }
  \label{fig:trader_control_MFG_t0}
\end{minipage}%
\hspace{0.2cm}\begin{minipage}{.45\textwidth} 
  \centering 
  \includegraphics[width=.9\linewidth]{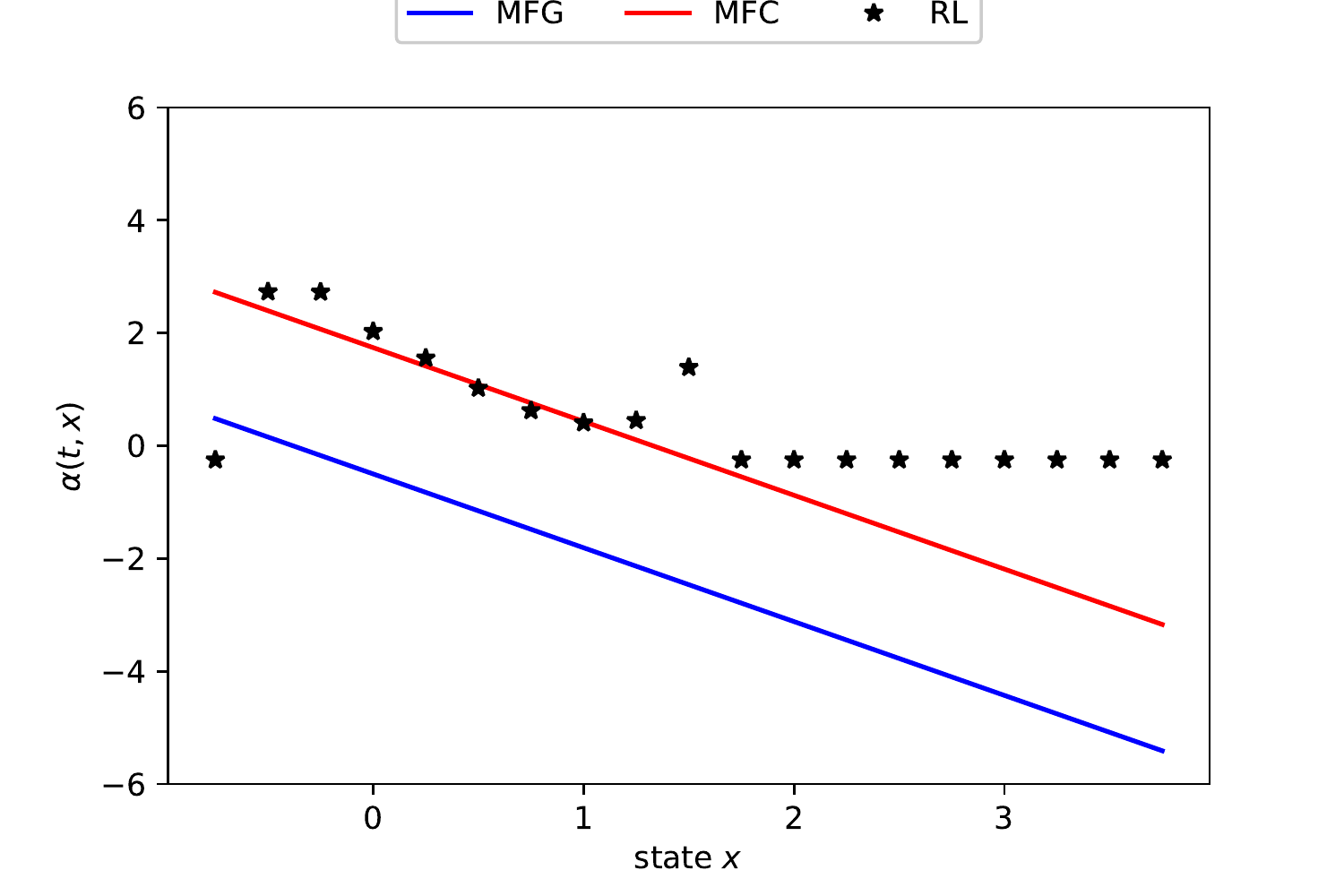}
  \caption{Learned Controls for MFC at time $0$. 
  }
  \label{fig:trader_control_MFC_t0}
\end{minipage}%
\end{figure}

\begin{figure}[H]
\centering
\begin{minipage}{.45\textwidth} 
  \centering 
  \includegraphics[width=.9\linewidth]{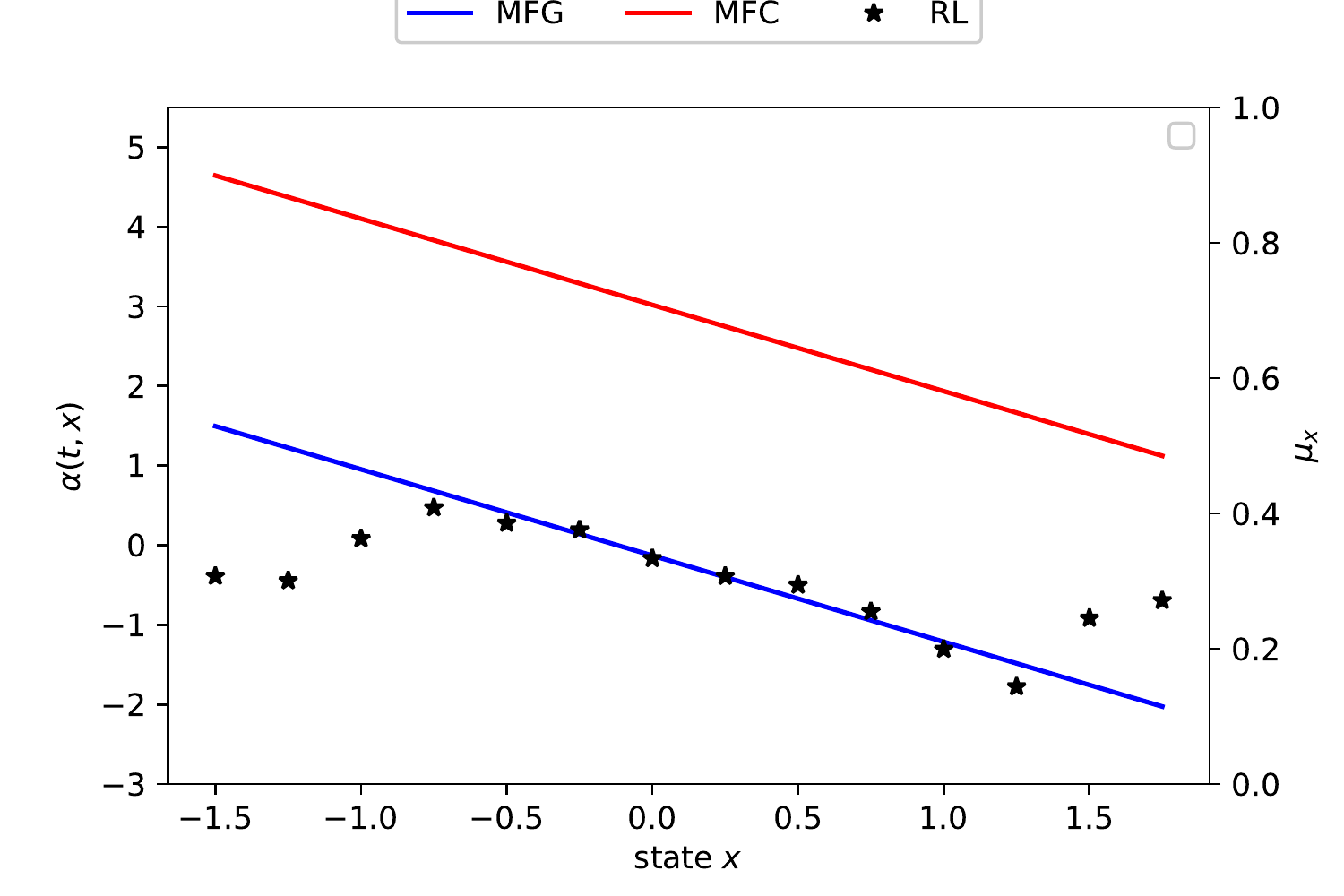}
  \caption{Learned Controls for MFG at time $7/16$.
  }
  \label{fig:trader_control_MFG_t50}
\end{minipage}%
\hspace{0.2cm}\begin{minipage}{.45\textwidth} 
  \centering 
  \includegraphics[width=.9\linewidth]{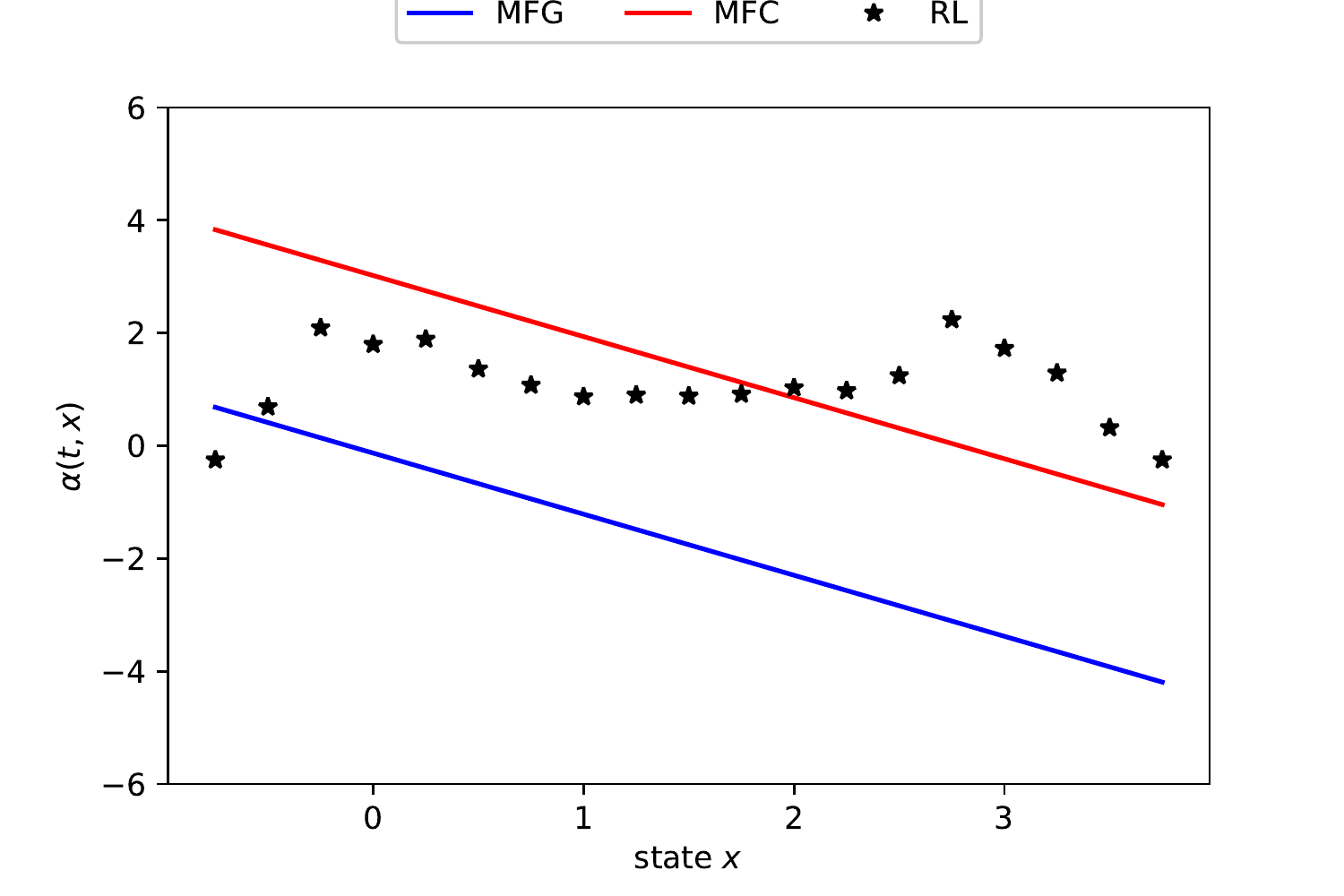}
  \caption{Learned Controls for MFC at time $7/16 $. 
  }
  \label{fig:trader_control_MFC_t50}
\end{minipage}%
\end{figure}

\begin{figure}[H]
\centering
\begin{minipage}{.45\textwidth} 
  \centering 
  \includegraphics[width=.9\linewidth]{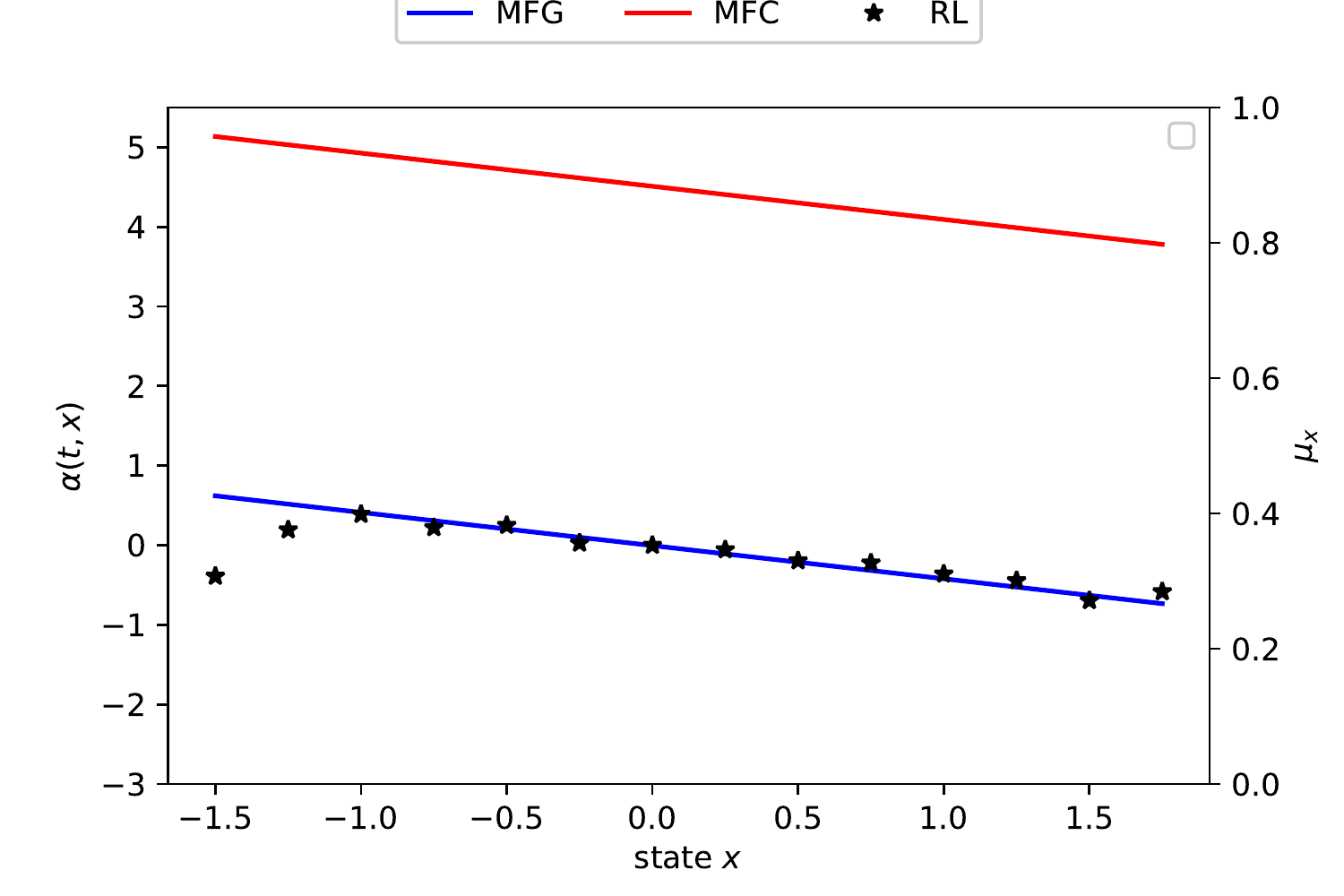}
  \caption{Learned Controls for MFG at time $15/16$.
  }
  \label{fig:trader_control_MFG_t100}
\end{minipage}%
\hspace{0.2cm}\begin{minipage}{.45\textwidth} 
  \centering 
  \includegraphics[width=.9\linewidth]{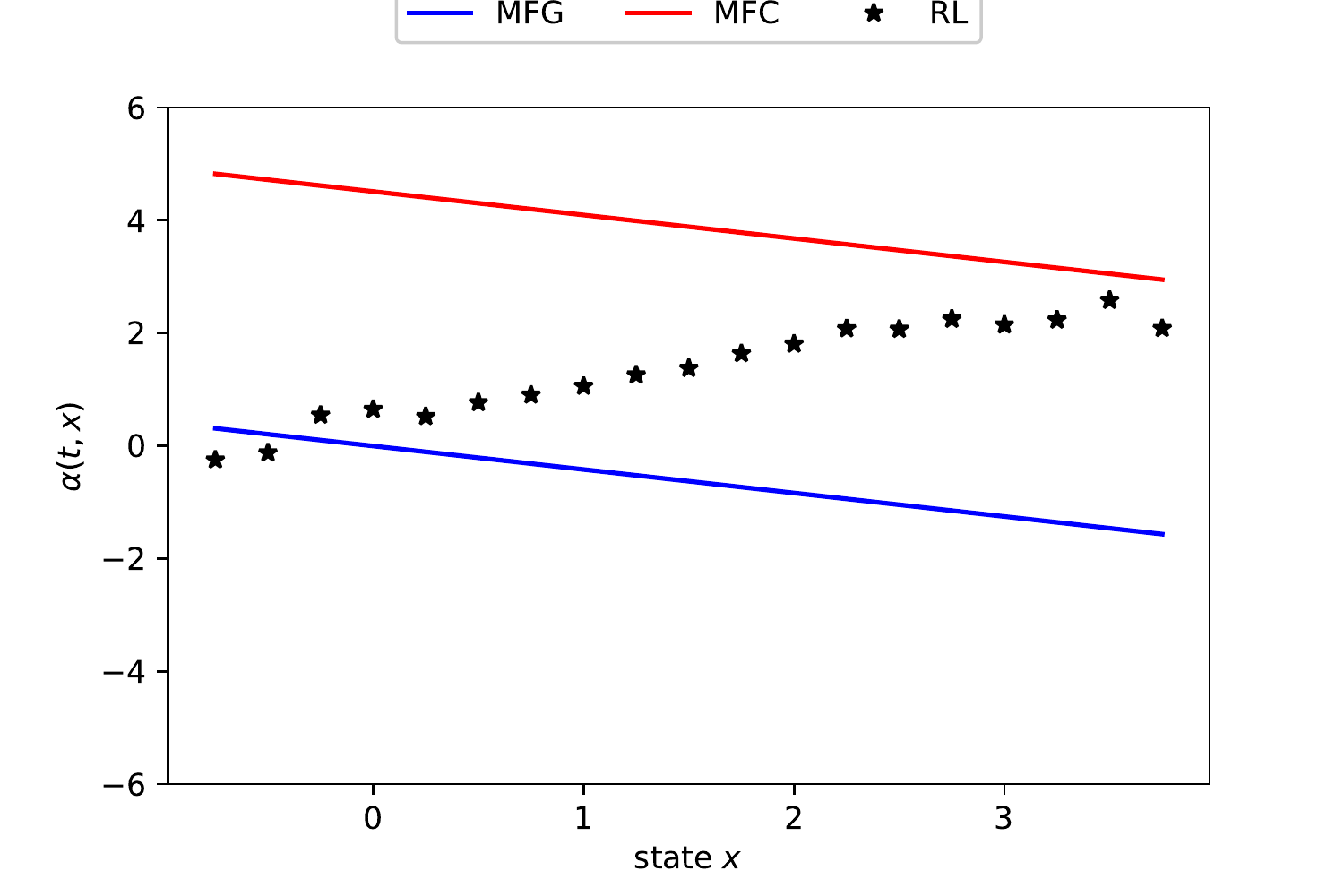}
  \caption{Learned Controls for MFC at time $15/16$. 
  }
  \label{fig:trader_control_MFC_t100}
\end{minipage}%
\end{figure}

\textbf{Figures \ref{fig:trader_control_MFG_t0}, \ref{fig:trader_control_MFC_t0},  \ref{fig:trader_control_MFG_t50},  \ref{fig:trader_control_MFG_t100}, \ref{fig:trader_control_MFC_t100}:  controls learned by the algorithm.} The controls learned by the U2-MF-QL-FH algorithm are compared with the theoretical solutions. Each plot corresponds to a different time point $t \in \{0,0.5,1\}$. The layout is the same applied for the mean field capital accumulation problem in Section \ref{sec:results_hara}. On the left, the choice $(\omega^Q,\omega^{\theta})=(0.55,0.85)$ produces the approximation of the solution of the MFG. On the right, the values of  the parameters  $(\omega^Q,\omega^{\theta})=(0.65,0.15)$ lets the algorithm to approach the solution of the MFC problem. The accuracy of the approximation is better at initial times and degrades towards the final horizon showing an higher complexity of the tuning of the algorithm to this problem. The results presented in the Figures are averaged over 10 runs.

\subsubsection{Learning of the mean field}
\begin{figure}[H]
\centering
\begin{minipage}{.45\textwidth} 
  \centering 
  \includegraphics[width=.9\linewidth]{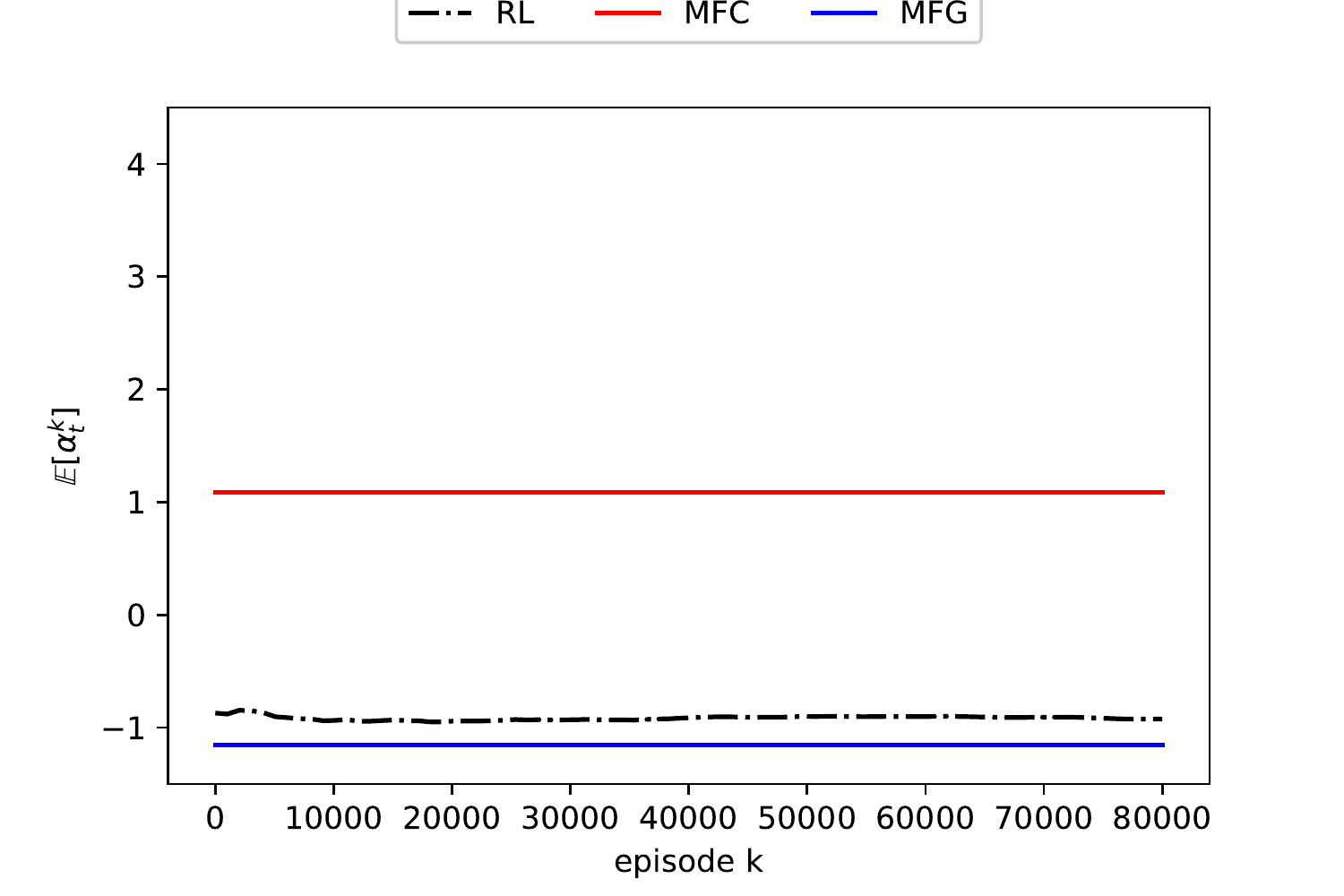}
  \caption{Learned control's mean for MFG at time $0$.
  }
  \label{fig:trader_theta_MFG_t0}
\end{minipage}%
\hspace{0.2cm}\begin{minipage}{.45\textwidth} 
  \centering 
  \includegraphics[width=.9\linewidth]{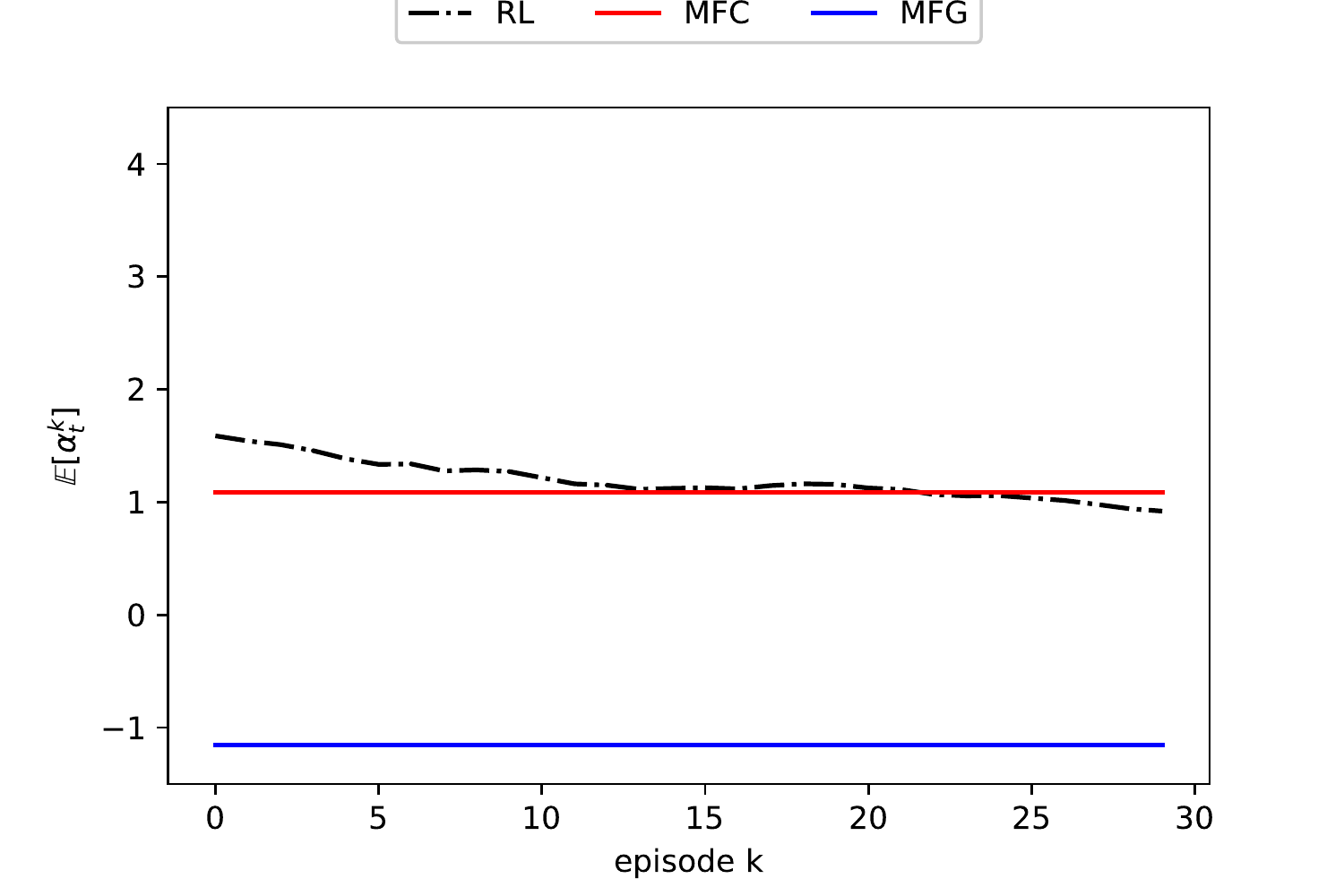}
  \caption{Learned control's mean for MFC at time $0$. 
  }
  \label{fig:trader_theta_MFC_t0}
\end{minipage}%
\end{figure}

\begin{figure}[H]
\centering
\begin{minipage}{.45\textwidth} 
  \centering 
  \includegraphics[width=.9\linewidth]{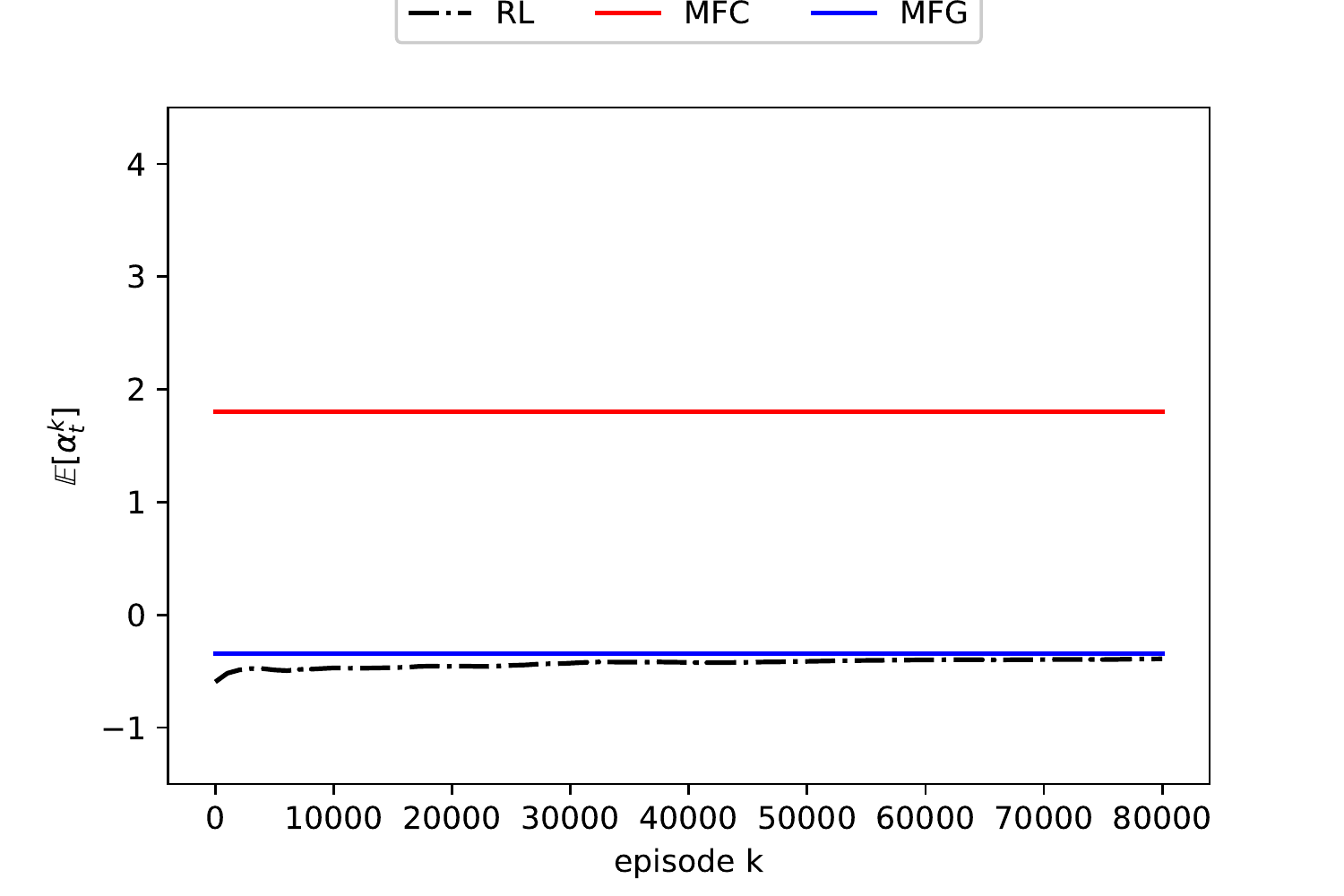}
  \caption{Learned control's mean for MFG at time $7/16$.
  }
  \label{fig:trader_theta_MFG_t50}
\end{minipage}%
\hspace{0.2cm}\begin{minipage}{.45\textwidth} 
  \centering 
  \includegraphics[width=.9\linewidth]{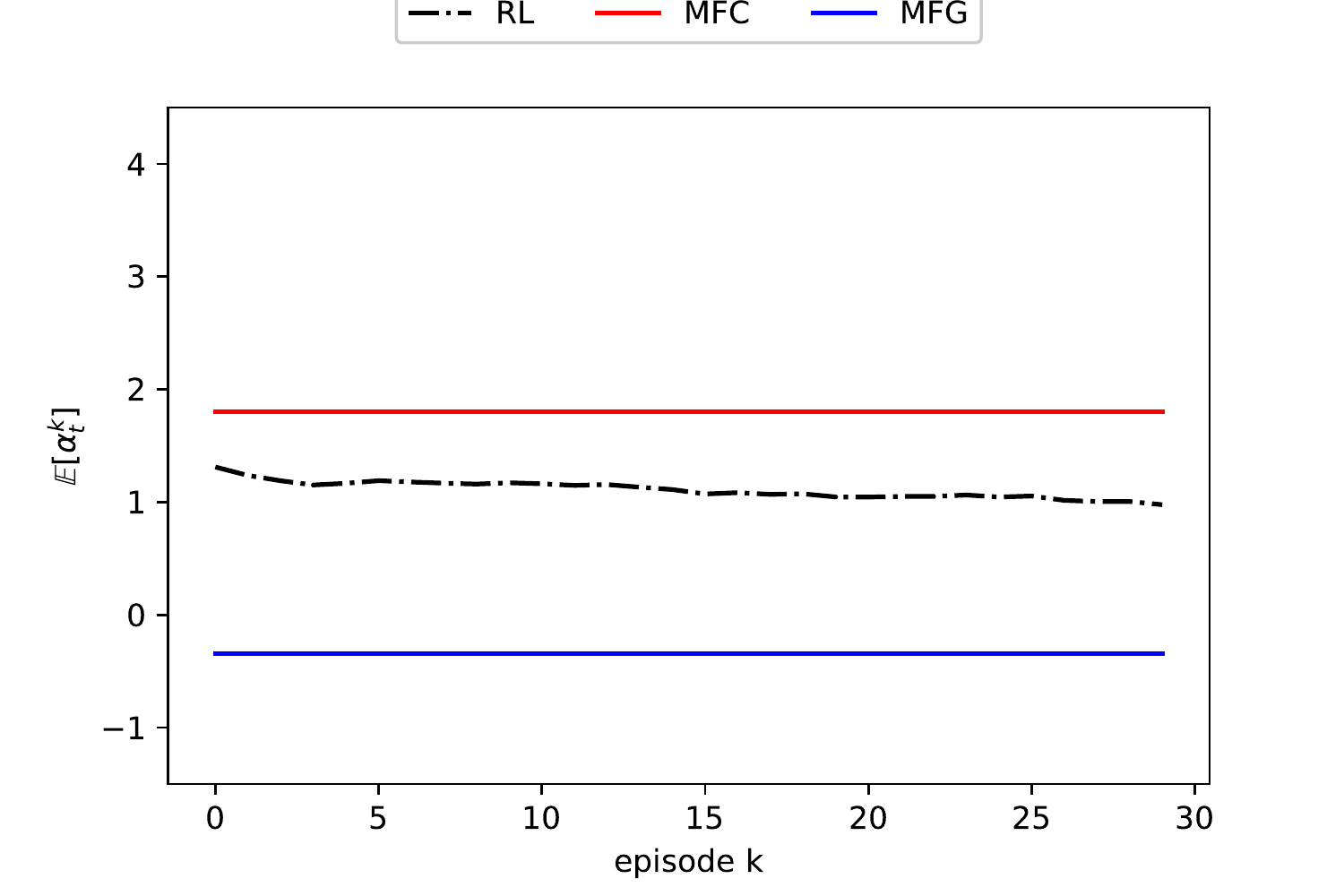}
  \caption{Learned control's mean for MFC at time $7/16 $. 
  }
  \label{fig:trader_theta_MFC_t50}
\end{minipage}%
\end{figure}

\begin{figure}[H]
\centering
\begin{minipage}{.45\textwidth} 
  \centering 
  \includegraphics[width=.9\linewidth]{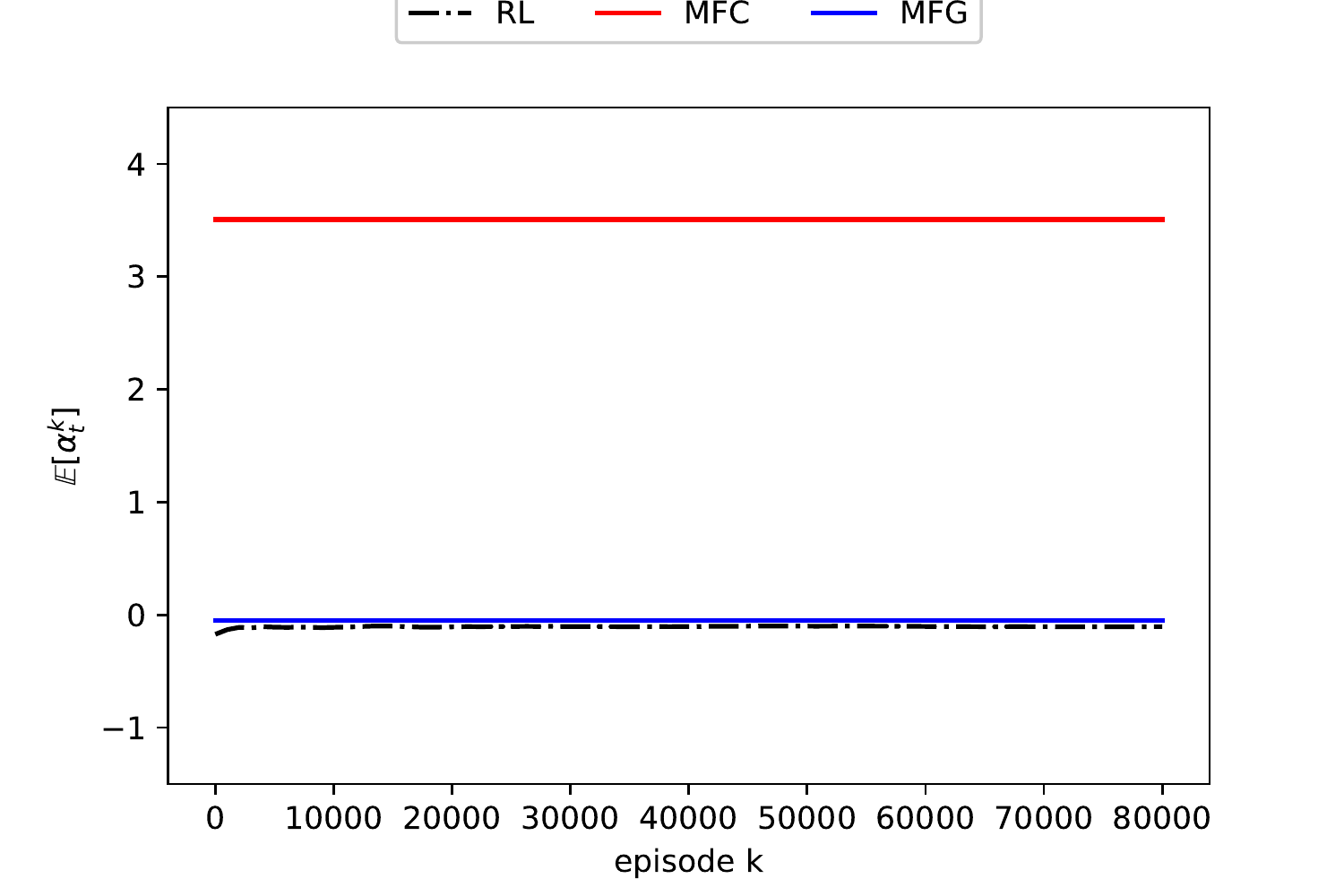}
  \caption{Learned control's mean for MFG at time $15/16$.
  }
  \label{fig:trader_theta_MFG_t100}
\end{minipage}%
\hspace{0.2cm}\begin{minipage}{.45\textwidth} 
  \centering 
  \includegraphics[width=.9\linewidth]{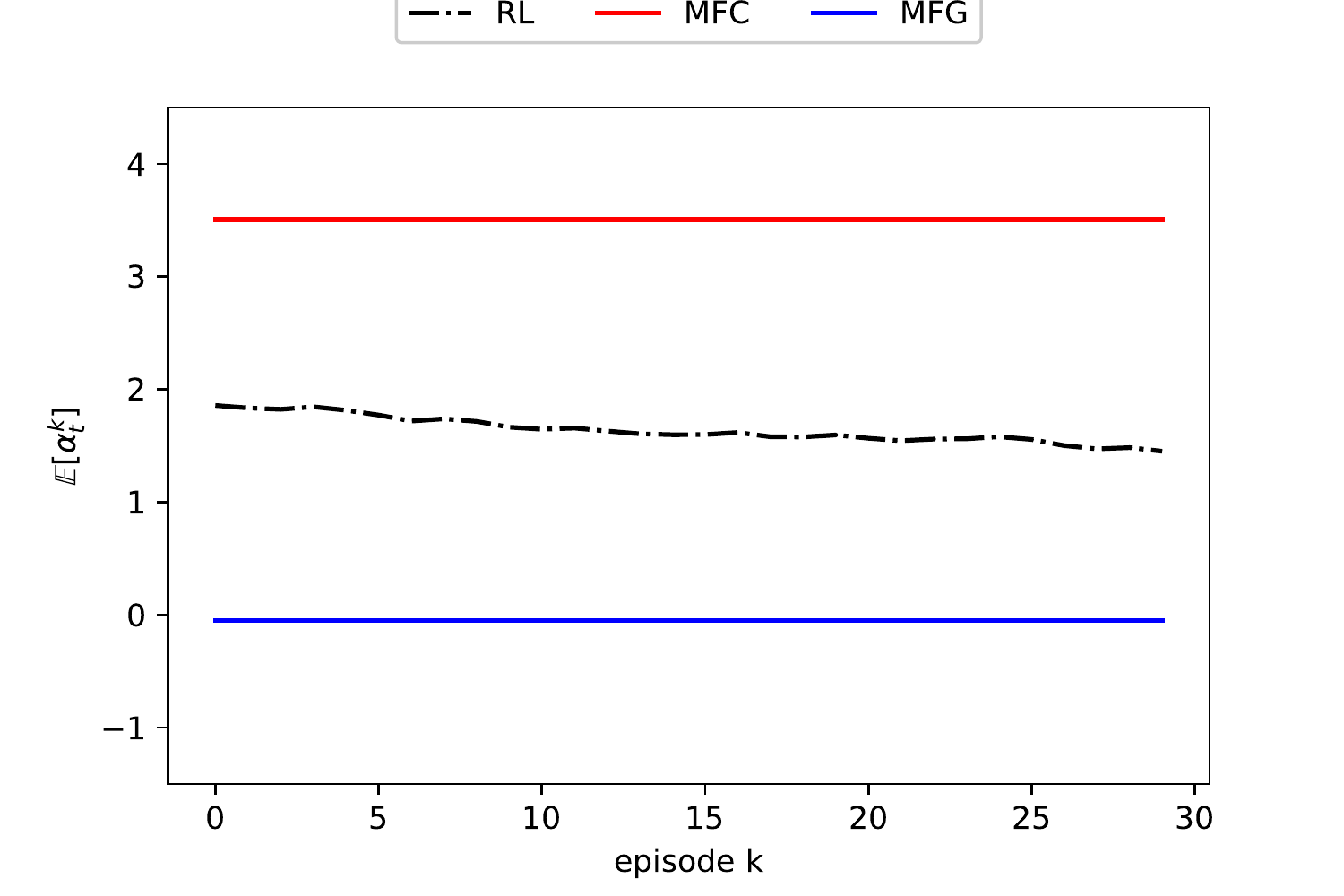}
  \caption{Learned control's mean for MFC at time $15/16$. 
  }
  \label{fig:trader_theta_MFC_t100}
\end{minipage}%
\end{figure}

\textbf{Figures \ref{fig:trader_theta_MFG_t0}, \ref{fig:trader_theta_MFC_t0},  \ref{fig:trader_theta_MFG_t50}, \ref{fig:trader_theta_MFC_t50},  \ref{fig:trader_theta_MFG_t100}, \ref{fig:trader_theta_MFC_t100}:  $\E[\theta_t]$ learned by the algorithm.} 
The estimation of the first moment of the distribution of the controls evolves with respect to the number of learning episodes. Each plot corresponds to a different time point $t \in \{0,0.5,1\}$. The layout is the same described in Section \ref{sec:results_hara}. On the left, the solution of the MFG is obtained choosing $(\omega^Q,\omega^{\theta})=(0.55,0.85)$. On the right, the MFC solution is approached by the set of parameters  $(\omega^Q,\omega^{\theta})=(0.65,0.15)$. The results presented in the Figures are averaged over 10 runs.

\section{Conclusion}
\label{sec:conclusion}

In this work, we have presented a reinforcement learning algorithm which can be used to approximate solutions of mean field game or mean field control problems in the case of interaction through the distribution of controls  and in finite horizon. The method unifies the two problems through a two timescale perspective. We have illustrated the algorithm with two examples: an optimal investment problem with HARA utility function, and 
an optimal liquidation problem. 

The main ingredients of the algorithm are the learning rates for the Q-matrix and for the distribution of controls. Their relative decay with respect to the number of episodes is the key quantity to stir the algorithm towards learning the optimal controls for MFG or MFC problem. Roughly speaking, updating the Q-matrix faster (resp. slower) than the distribution of controls leads to the MFG (resp. MFC) solution. 
Convergence  follows by applying Borkar's results as shown in \citep{angiuli2020unified} in the case of infinite horizon problems.
Choosing these rates in an optimal way remains the main challenge in specific applications. In particular, we expect that allowing these rates to depend on the time steps could lead to improved results. This aspect is left for future investigations.

The algorithm presented here is the context of finite space via the Q-matrix even though the proposed examples are originally in continuous space and then discretized. Dealing directly with a continuous space is the topic of the ongoing work on deep reinforcement learning for mean filed problems \citep{AngiuliHU-2021-deep_mfrl}.

The area of reinforcement learning for mean field problems is extremely rich with a huge potential for applications in various disciplines. It is in its infancy, and we hope that the results and explanations presented here will be helpful to newcomers interested in  this direction of research.

\bibliographystyle{apalike}
\bibliography{references}

\end{document}